\documentclass[oo,nonblindrev]{informs3}

\usepackage{amssymb}
\usepackage{amsmath}
\usepackage{mathtools}
\usepackage{bm}

\usepackage{algorithm}
\usepackage[noend]{algpseudocode}

\usepackage{enumerate}
\usepackage{enumitem}

\usepackage[table]{xcolor}
\usepackage{graphicx} 
\usepackage{adjustbox}
\usepackage{lscape} 
\usepackage{subcaption}
\usepackage{comment}

\usepackage{multirow}
\usepackage{booktabs}

\usepackage{appendix}
\usepackage{float}
\usepackage{siunitx}

\usepackage{url}
\usepackage{hyperref}
\hypersetup{colorlinks=true,linkcolor=black,citecolor=black,filecolor=black,urlcolor=black,pdfpagelayout=SinglePage}

\usepackage{tikz}
\usetikzlibrary{plotmarks}
\usetikzlibrary{arrows}

\definecolor{kellygreen}{rgb}{0.3,0.73,0.09}
\definecolor{lightgreen}{rgb}{0.56,0.93,0.56}
\definecolor{lightyellow}{rgb}{0.94,0.91,0.57}
\definecolor{lightorange}{rgb}{1.0,0.85,0.69}
\definecolor{lightsalmon}{rgb}{1.0,0.63,0.48}
\definecolor{lightcoral}{rgb}{0.94,0.5,0.5}

\usepackage{endnotes}
\let\footnote=\endnote

\usepackage{natbib}
\bibpunct[, ]{(}{)}{,}{a}{}{,}%
\def\bibfont{\small}%
\TheoremsNumberedThrough     
\ECRepeatTheorems
\EquationsNumberedThrough    

\usepackage{graphicx}
\makeatletter
\newcommand*\bigcdot{\mathpalette\bigcdot@{.5}}
\newcommand*\bigcdot@[2]{\mathbin{\vcenter{\hbox{\scalebox{#2}{$\m@th#1\bullet$}}}}}
\makeatother

\usepackage{enumerate}
\usepackage{pgfplots}
\usepackage{arydshln}
\usepackage[official]{eurosym}
\usepackage{lineno}

\usepackage{amsmath}
\newcommand{\decision}{p}
\newcommand{\product}{W}
\newcommand{\denuminator}{f}
\newcommand{\ahmadreza}[1]{{\color{black}#1}}
\newcommand{\probInd}{\ell}
\newcommand{\probNum}{L}

\OneAndAHalfSpacedXII

\TITLE{An exact algorithm for the static pricing problem under discrete mixed logit demand}

\ARTICLEAUTHORS{
\AUTHOR{Ahmadreza Marandi}
\AFF{Eindhoven University of Technology, Department of Industrial Engineering \& Innovation Sciences, The Netherlands.}
\AFF{Eindhoven University of Technology, Eindhoven Artificial Intelligence Systems Institute, The Netherlands, \EMAIL{a.marandi@tue.nl}.}
\AUTHOR{Virginie Lurkin}
\AFF{University of Lausanne, HEC Lausanne, Faculty of Business and Economics, Switzerland, \EMAIL{virginie.lurkin@unil.ch}.}
}

\begin{document}

	\RUNAUTHOR{A. Marandi, V. Lurkin}

	\RUNTITLE{An exact algorithm for the static pricing problem under discrete mixed logit demand}
\ABSTRACT{
Price differentiation is a common strategy in many markets. In this paper, we study a static multiproduct price optimization problem with demand given by a discrete mixed multinomial logit model. By considering a mixed logit model that includes \ahmadreza{customer} specific variables and \ahmadreza{parameters} in the utility specification, our pricing problem reflects well the discrete choice models used in practice. To solve this pricing problem, we design an efficient iterative optimization algorithm that asymptotically converges to the optimal solution. To this end, a \textit{linear optimization} (LO) problem is formulated, based on the trust-region approach, to find a ``good" feasible solution and approximate the problem from below. A \ahmadreza{convex optimization} problem is designed using a \ahmadreza{convexification technique} to approximate the optimization problem from above. Then, using a branching method we tighten the optimality gap. The effectiveness of our algorithm is \ahmadreza{illustrated} on several cases, and \ahmadreza{compared} against solvers and existing \ahmadreza{state-of-the-art methods} in the literature. 
}


\maketitle
\section{Introduction}
Offering different products at different prices to different customers is a common practice in many markets, including transportation, retail\ahmadreza{,} and entertainment, among others. Classic examples include discount prices for senior\ahmadreza{s} and children, business-, first-, and economy-class flight tickets, or first- and second-class railway tickets. With product and price differentiation, operators and retailers are able to get higher revenues by adapting their \ahmadreza{prices} based on the price sensitivity of their customers. Basically, higher fares are offered to the ones who are willing to pay more.

Inferring customers’ willingness to pay (WTP) is a long-standing practice in applied economics \citep{hensher2005applied}. Discrete-choice modeling (DCM) has established itself as an important and widely-used methodology for extracting valuations such as willingness to pay \citep{hess2018revisiting}. Researchers have used these disaggregate demand models for more than 40 years, from the pioneering work of \cite{mcfadden1974frontiers} to more recent studies on WTP for self-driving vehicles \citep{daziano2017consumers} or willingness to travel with green modes in the context of shared mobility \citep{li2019integrated}. 

Formulating pricing policies based on such disaggregate demand representations allows to better account for the heterogeneity of the population of interest, where different customers have different tastes and preferences. Even more importantly, it better reflects the supply-demand interactions by capturing the tradeoff between the operator\ahmadreza{'s} objective of maximizing the expected revenue and the customer objective of maximizing the expected utility \citep{sumida2019revenue}. 

Despite a more comprehensive representation, including discrete choice models within pricing problems increases the computational complexity because the choice probabilities are nonlinear. As a result\ahmadreza{,} the expected revenue is highly nonlinear in the prices of the products\ahmadreza{,} and customary used nonlinear algorithms may get terminated at a local optimum.

Due to the importance of the problem, the operations research and management science communities put remarkable efforts \ahmadreza{into} analyzing it. \cite{hanson1996} pioneer this research by showing that the expected revenue function is not concave in prices, even for the simple \textit{multinomial logit} (MNL) model. Subsequent authors have demonstrated that, under uniform price sensitivities across all products, the expected revenue function is concave in the choice probability vector \citep{song2007demand,dong2009dynamic,zhang2013assessing}. \cite{li2011pricing} show that this concavity result also holds under asymmetric price-sensitivities, not only for the MNL model, but also for the \textit{nested logit} (NL) model that generalizes the MNL model by grouping product alternatives into different nests based on their degree of substitution \citep{mcfadden78}. 

 Parallel to these works, several authors have shown that under restrictive conditions on the degree of asymmetry in the price sensitivity parameters, unique price solutions exist for some logit models. This has been shown for the MNL model (e.g., \cite{aydin2000product,hopp2005product,maddah2007joint,aydin2008joint,akccay2010joint}), the NL model (e.g., \cite{aydin2000product,hopp2005product,maddah2007joint,aydin2008joint,akccay2010joint,gallego2014multiproduct,huh2015pricing}), the \textit{paired combinatorial logit} (PCL) model \citep{li2017optimal} and lately generalized to any \textit{generalized extreme value} (GEV) model \citep{zhang2018multiproduct}. In this stream of research, a first-order condition is generally used to find optimal prices. It is \ahmadreza{worth noting} that in some of these studies and additional recent ones, (1) pricing decisions are optimized jointly with other decisions such as assortment or scheduling decisions (e.g., \cite{du2016optimal,jalali2019quality,bertsimas2020joint}), (2) decisions of multiple firms (or players) are studied, mainly from a non-cooperative game theory perspective (e.g., \cite{li2011pricing,aksoy2013price,gallego2014multiproduct,bortolomiol2021simulation}), but also more recently, from a cooperative game theory perspective \citep{schlicher2021stable}.
 
Only a few papers (\cite{gilbert2014mixed}, \cite{li2019product}, and \cite{GeerBranch}) consider pricing problems under a \textit{mixed logit} (ML) model, a choice model that better accommodates customer heterogeneity by allowing some parameter\ahmadreza{s} to vary across \ahmadreza{customer}s. As \ahmadreza{shown} by \cite{mcfadden2000mixed}, under mild regularity conditions, the ML model can approximate choice probabilities of any discrete choice model derived from random utility maximization (RUM) assumption, making it a popular choice model. However, as explained by \cite{li2019product}, the expected revenue function under the mixed logit model is not well-behaved\ahmadreza{,} and the concavity property with respect to the choice probabilities breaks down, even for entirely symmetric price sensitivities across products and segments. Accordingly, the theoretical results as well as the solution methods developed for other logit models do not apply to the pricing problem with demand characterized by a mixed logit model. 

In \cite{gilbert2014mixed}, the authors consider a ML demand model within a revenue-maximizing network pricing problem whose objective is to improve the performance of a congested network through the selection of appropriate tolls. The price sensitivity parameter is distributed across the population according to a continuous random variable. The authors rely on a tractable approximation of the ML pricing problem. Their approach consists of two phases. They first solve optimally a mixed integer program that approximates the original problem by assuming a simpler distribution for \ahmadreza{the} price sensitivity parameter. The optimal solution of this program is then considered as the starting solution of an ascent algorithm that is used to solve a differentiable optimization problem that better approximates the original ML pricing problem.  
 
In \cite{li2019product}, the authors assume that the market is divided into a finite number of market segments, with product demand in each segment governed by the multinomial logit model. The problem under investigation is therefore a price optimization problem with demand given by a \textit{discrete} ML model. To solve this problem, the authors propose two concave maximization problems that work as lower and upper bounds for the objective value of the revenue function, under some conditions. Then, they propose an algorithm that converges to a local optimum. 

In \cite{GeerBranch}, the authors also assume that the market is divided into a finite number of market segments and \ahmadreza{that} customers’ intrinsic product valuations
are both product and segment dependent. However, the customers’ price sensitivity parameters are product dependent only,  and thus identical for all customer segments. This assumption was critical for them to develop a scalable algorithm that quickly converges to an optimal price of products. \color{black}

Our pricing problem also considers the discrete setting of the ML model and therefore includes \ahmadreza{customer} specific variables in the utility specifications.  However, unlike former contributions that either consider local optimality \citep{hanson1996} or impose restrictive conditions on price sensitivity parameters (\cite{li2019product} and \cite{GeerBranch}) to have global optimality, our pricing problem can handle customers’ price sensitivity parameters that differ among products and customer segments. Considering that customers have different willingness to pay (WTP) is more realistic, but it comes at the price of additional computational complexity. 

To solve our pricing problem,\color{black} we follow an approach similar to \cite{GeerBranch} by designing an exact iterative optimization algorithm that asymptotically converges to \ahmadreza{a} global optimal solution. More specifically, we develop a method to find ``good'' solutions, design an approach to check the quality of the obtained solution, and branch to make sure the solution is optimal. To this end, a \textit{linear optimization} (LO) problem is formulated, based on the trust-region approach, to find a ``good" feasible solution and approximate the problem from below. A \ahmadreza{convex optimization} problem is designed using \ahmadreza{a convexification technique} as well as the McCormick relaxation \citep{mccormick1976computability} to approximate the optimization problem from above.  Then, we develop a branching method to tighten the optimality gap and show that the algorithm converges to the optimal solution asymptotically. The effectiveness of this algorithm is \ahmadreza{illustrated} on several instances, including a parking services pricing case for which the demand model comes from a published, non-trivial parking choice model.

The remaining sections are organized as follows. Section 2 further defines the problem under consideration. Section 3 presents our global algorithm, while Section 4 shows the results of our numerical experiments. The final section concludes our paper.

\section{Problem description}
In this paper, we are interested in solving a static multi-product pricing problem under a discrete ML model. Static pricing involves the simultaneous pricing of multiple products, where a fixed price is set for each product \citep{soon2011review}. In our setting, we assume that a single seller must decide at what price to offer each product from a finite set of alternatives (also known as product assortment). On the demand side, we assume that customers choose among the products according to a consumer choice model. The demand for each product is thus the result of the \ahmadreza{customer} purchase choice of $N$ customers. The purchase choice is captured by a discrete choice model, that predicts the customer choice from a finite set of discrete alternatives \citep{ben2003discrete}.

Let $\mathcal{N}$ represent the set of $N$ customers and  let $\mathcal{\ahmadreza{C}}$ indicate the set of $C$ available alternatives, among which $I$ alternatives are offered by the seller. We denote by $\mathcal{I}$ the offered alternatives. So, we assume $\mathcal{I}\subseteq \mathcal{C}$, and  $\mathcal{C}\setminus \mathcal{I}$ is the set of other alternatives not offered by the seller, such as competitive products or the standard non-purchase alternative. \color{black} \ahmadreza{For each customer $n \in \mathcal{N}$ and  alternative $i \in \mathcal{C}$, the} utility function $U_{in}$ is a function of the socio-economic characteristics of the \ahmadreza{customer}s and/or the attributes of the alternatives. According to Random Utility Maximization (RUM) theory \citep{manski1977structure}, $U_{in}$ can be decomposed into a systematic component $V_{in}(\beta)$, which includes all attributes observed by the decision maker\ahmadreza{, $\beta$,} and a random term $\varepsilon_{in}$, which captures the uncertainties caused by unobserved attributes and unobserved taste variations: 
\begin{align}\label{Eq:utility}
U_{in} &= V_{in}(\beta) + \varepsilon_{in}= \beta^p_{in} p_{i} + q_{in}(\beta^q) + \varepsilon_{in},
\end{align}
where $p_{i}$ is the price of the alternative $i$, $\beta^p_{in}$ is the willingness to pay of customer $n$ for alternative $i$, and $q_{in}(\beta^q)$ is the exogenous part of the utility, obtained by adding all observed product attributes other than price, weighted based on customers' preferences.  Note that for the alternatives that are not offered by the seller, the price is assumed to be fixed and given, i.e., $p_{i} =\bar{\decision}_i, \forall i \in \mathcal{C}\setminus \mathcal{I}$. \color{black}

The resulting discrete choice model is\ahmadreza{,} therefore\ahmadreza{,} naturally probabilistic. The probability that customer $n$ chooses alternative $i$ is defined as 
$$P_{in} = \Pr\left[V_{in}(\beta) + \varepsilon_{in} = \max_{j \in \mathcal{I}} \left\{V_{jn}(\beta) + \varepsilon_{jn}\right\}\right].$$

The optimal expected revenues \ahmadreza{of the seller}, obtained from the sales of \ahmadreza{all product within the set $\mathcal{I}$}, is then given by: 
\begin{equation}\label{Eq:choicebasedpricing}
\begin{aligned}
\max_{p\in\mathbb{R}^{\ahmadreza{C}}}& \enskip \sum_{i\in\mathcal{I}} \sum_{n \in \mathcal{N}} p_i P_{in}, & \\
\mbox{s.t.}&\enskip P_{in} = \Pr\left[V_{in}(\beta) + \varepsilon_{in} \geq V_{jn}(\beta) + \varepsilon_{jn},\enskip \forall j \in \mathcal{\ahmadreza{C}}\right], & \forall i \in \mathcal{\ahmadreza{C}}, \forall n \in \mathcal{N}, \\
\mbox{}&\enskip V_{in}(\beta) = \beta^p_{in} p_{i} + q_{in}(\beta^q),&\forall i \in \mathcal{\ahmadreza{C}}, \forall n \in \mathcal{N}, \\
\mbox{}&\enskip  \ahmadreza{p_i=  \bar{\decision}_i}, &\forall i \in \mathcal{\ahmadreza{C \setminus I}}, \\
\mbox{}&\enskip 0\leq p_i\leq \ahmadreza{{\decision}^u_i}, &\forall i \in \mathcal{\ahmadreza{I}},
\end{aligned}
\end{equation}
where $\ahmadreza{\decision^u}\in\mathbb{R}^I$ is a vector containing upper bounds on the prices of products of the seller \ahmadreza{and $\bar{\decision}_i$ is the given price of the exogenous alternative $i \in\mathcal{C}\setminus \mathcal{I}$}. 
The most commonly used discrete choice models, the multinomial logit (MNL) model, is built upon the assumption of independent and identically extreme value distributed error terms \citep{manski1977structure}, that is $\varepsilon_{in} \overset{\text{i.i.d.}}{\sim} EV(0, 1)$. Under this assumption, the probability for customer $n$ to select choice alternative $i$ is given by
\begin{equation}
\label{eq:MNLproba}
P_{in} = \frac{e^{V_{in}(\beta)}}{\sum_{j\in \mathcal{\ahmadreza{C}}} e^{V_{jn}(\beta)}}. 
\end{equation}

Mixed logit probabilities are the weighted sum of these standard logit probabilities over a density of parameters \citep{train_2003}. The choice probabilities can then be expressed as

\begin{equation}
\label{eq:mixedproba}
P_{in} = \int \frac{e^{V_{in}(\beta)}}{\sum_{j\in \mathcal{\ahmadreza{C}}} e^{V_{jn}(\beta)}} d\nu_{\beta},
\end{equation}
where $\nu_{\beta}$ is a multivariate probability measure. 

In this paper, we assume that the probability measure $\nu_{\beta}$ is discrete and can be dependent on \ahmadreza{customer}s. In other words, \ahmadreza{$\beta$} can take only $\probNum$ distinct values \ahmadreza{$b_1,\enskip b_2, \; ...,\; b_{\probNum}$}, probability of which differs per \ahmadreza{customer}. So, we have the following logit choice probability:
\begin{equation}
\label{eq:mixedprobadiscrete}
P_{in} = \sum_{\probInd=1}^{\probNum} w_{\probInd n} \frac{e^{V_{in}(\ahmadreza{b_{\probInd}})}}{\sum_{j\in \mathcal{\ahmadreza{C}}} e^{V_{jn}(\ahmadreza{b_{\probInd}})}},
\end{equation}
where $w_{\probInd n}$ is the probability that $\beta^p = b^p_{\probInd}$ for \ahmadreza{customer} $n$.   


We are thus interested in solving the following nonlinear maximization problem: 

\begin{equation}\label{Eq:ourproblem}
\begin{aligned}
\max_{p\in\mathbb{R}^{\ahmadreza{C}}}& \enskip \sum_{i\in\mathcal{I}} \sum_{n \in \mathcal{N}} p_i P_{in}, & \\
\mbox{s.t.}&\enskip P_{in} = \sum_{\probInd=1}^{\probNum} w_{\probInd n} \frac{e^{V_{in}(\ahmadreza{b_{\probInd}})}}{\sum_{j\in \mathcal{\ahmadreza{C}}} e^{V_{jn}(\ahmadreza{b_{\probInd}})}}, &\forall i \in \mathcal{\ahmadreza{C}}, \forall n \in \mathcal{N}, \\
\mbox{}&\enskip V_{in}(\beta) = \beta^p_{in} p_{i} + q_{in}(\beta^q),&\forall i \in \mathcal{\ahmadreza{C}}, \forall n \in \mathcal{N}, \\
\mbox{}&\enskip  \ahmadreza{p_i=  \bar{\decision}_i}, &\forall i \in \mathcal{\ahmadreza{C \setminus I}}, \\
\mbox{}&\enskip 0\leq p_i\leq \ahmadreza{{\decision}^u}, &\forall i \in \mathcal{\ahmadreza{I}},
\end{aligned}
\end{equation}


It is worth noting that in \cite{li2019product} and \cite{GeerBranch}, the closest related works, the probability measure $\nu_{\beta}$ is also assumed to be discrete. In \cite{li2019product}, there is no \ahmadreza{customer} specific variables included in their utility specifications. As a result, the choice probability is the same for all \ahmadreza{customer}s (i.e., $P_{in} = P_i, ~\forall n \in N$) and the objective function becomes $\max_{p\in\mathbb{R}^I}  \sum_{i\in\mathcal{I}}  p_i P_{i}$. In \cite{GeerBranch}, only the exogenous part of the utility is linked to the \ahmadreza{customer}s. So, unlike \cite{GeerBranch}, our pricing problem handles customer heterogeneity also in terms of price sensitivity parameters.



\section{Methodology}
\label{methodo}
In this section, we introduce a new efficient optimization algorithm for solving Problem \ahmadreza{\eqref{Eq:ourproblem}}. The proposed algorithm is a global optimizer, meaning that it asymptotically converges to the optimal solution. This is done by designing a method to find a ``good'' feasible solution, which provide a lower bound, as well as a method to check the quality of the obtained solution, which provides an upper bound. 

Let us  reformulate the optimization problem \eqref{Eq:ourproblem} to
\begin{equation}\label{Eq:generalized problem}
\begin{aligned}
opt=\max_{\decision\in\mathbb{R}^\ahmadreza{C}}& \enskip f(\decision)\\
\mbox{s.t.}&\enskip A\decision\geq b,\enskip \decision\geq 0,
\end{aligned}
\end{equation}
where $f\ahmadreza{:\mathbb{R}^I\rightarrow\mathbb{R}}$ is $\sum_{i\in\mathcal{I}}\sum_{n\in\mathcal{N}}\sum_{\probInd=1}^\probNum \frac{{w_{\probInd n}}\decision_i}{\denuminator_{i{\probInd n}}(\decision)}$, $\denuminator_{i{\probInd n}}:\mathbb{R}^I\rightarrow\mathbb{R}_\ahmadreza{+}$ is ${\sum_{j\in \ahmadreza{\mathcal{C}}} e^{V_{jn}(\ahmadreza{b_{\probInd}})-V_{in}(\ahmadreza{b_{\probInd}})}},$  for any $i\in\mathcal{I},$ $n\in\mathcal{N},$ and $\probInd=1,...,\probNum,$ \ahmadreza{$\mathbb{R}_+$ is the set of positive real numbers,} and the feasible region is a polytope defined by  the matrix $A\in\mathbb{R}^{m\times \ahmadreza{C}}$ and the vector $b\in\mathbb{R}^m$. 


\subsection{Designing a method to construct lower bounds}\label{Sec:lower bound}
To construct the lower bounds, we use a trust-region method \citep{conn2000trust}, where solutions are obtained iteratively in the neighborhood of the previous feasible solution. A typical way of finding a better solution is by approximating the objective function with a quadratic function and solving the following optimization problem in the $k^{\mbox{th}}$ iteration:
\begin{equation}\label{Eq:typical trust region}
\begin{aligned}
\max_{\decision\in\mathbb{R}^{\ahmadreza{C}}}&\enskip \frac{1}{2}\decision^TH_k\decision+g_k^T\decision\\
\mbox{s.t.}&\enskip \Vert \decision-\decision^k\Vert_2\leq r_k\\
&\enskip A\decision\geq b,\enskip \decision\geq 0,
\end{aligned}
\end{equation}
where $H_k\in\mathbb{R}^{\ahmadreza{C\times C}}$ is the Hessian matrix and $g_k\in\mathbb{R}^{\ahmadreza{C}}$ is the gradient vector of the objective function at the feasible solution  $\decision^k$ obtained in the $(k-1)^{\mbox{st}}$ iteration, $r_k$ is the radius of the neighborhood, and where $\Vert.\Vert_2$ is the Euclidean norm. The issue is that the objective function of \eqref{Eq:generalized problem} is neither convex nor concave, and hence \eqref{Eq:typical trust region} might be a nonconcave quadratic optimization problem, known to belong to the class of NP-hard problems \citep{pardalos1991quadratic}. To avoid this issue, we use the linear approximation of the objective function in each iteration and use the following optimization problem:
\begin{equation}\label{Eq:our trust region}
    \begin{aligned}
    \max_{\ahmadreza{p\in\mathbb{R}^C}}&\enskip g_k^T\decision\\
    \mbox{s.t.}&\enskip \Vert \decision-\decision^k\Vert_1 \leq r_k\\
&\enskip A\decision\geq b,\enskip \decision\geq 0,
    \end{aligned}
\end{equation}
where $\Vert.\Vert_1$ is the $\ell_1$-norm. 

Algorithm \ref{alg:trust region} provides the steps taken to find a ``good'' feasible solution using \eqref{Eq:our trust region}. As one can see, \eqref{Eq:our trust region} is a linear optimization problem; hence optimal solutions are in the boundary points of its feasible region. So, the algorithm starts with searching for a good solution in the boundary of the neighborhood of the initial solution with \ahmadreza{a} radius $1$. It continues the search unless it does not reach a point with improvement in the objective function. Then, the radius of the neighborhood gets decreased with the hope of finding a better solution (\ahmadreza{in the numerical results, we set the decreasing scale to 10; i.e., we multiply the radius with 0.1}). The algorithm gets terminated when the improvement in the last two iterations \ahmadreza{is} less than a given tolerance error $\theta$, hence a local optimum.  We emphasize that we chose \ahmadreza{the} trust region approach as it is known to have an extremely fast convergence rate to a local optimum \citep{higham1999trust}. However, any other method to efficiently find a ``good'' solution works.

\begin{algorithm}
\caption{Steps to obtain a ``good'' feasible solution using \eqref{Eq:our trust region}}\label{alg:trust region}
\begin{algorithmic}[1]
\State \ahmadreza{Input: $A, \; b,\; r^0, \; \theta$, and the gradient of $f(p);$}
\State   select a random feasible solution $\decision^0$
\State $f^1:=+\infty$, $r^0:=1$, $k=0,$
 \While{$|f^1-f(\decision^{k})|>\theta$, }
\State   find $\decision^{k+1}$ by solving \eqref{Eq:our trust region} with  radius $r^0$
  \State $\bar{\decision}^0\leftarrow \decision^k$, $\bar{\decision}^1\leftarrow \decision^{k+1}$, $f^0\leftarrow f(\bar{\decision}^{0})$, $f^1\leftarrow f(\bar{\decision}^{1})$
  \While{$f^1>f^0$}\label{algSteps:while improvement}
  \State$\bar{\decision}^0\leftarrow \bar{\decision}^1$, $f^0\leftarrow f^1$
  \State find $\bar{\decision}^1$ by solving \eqref{Eq:our trust region} with initial point $\bar{\decision}^0$ and radius $r^0$
  \State $f^1\leftarrow f(\bar{\decision}^1)$, $r^0\leftarrow 1$
  \EndWhile\label{algSteps:while improvement end}
  \State$r^0\leftarrow \frac{r^0}{10}$\label{algstep:shortening step size}
  \State $p^{k+1}\leftarrow \bar{\decision}^1$, increase $k$ by $1$
 \EndWhile
\State \textbf{return} $\decision^k$
\end{algorithmic}
\end{algorithm}

\subsection{Designing a method to construct upper bounds}\label{Sec:upper bound}
In this section, we explore the properties of the optimization problem \eqref{Eq:generalized problem} and use them to develop an overestimator to construct an upper bound on the objective value of the problem. To this end, we first reformulate \eqref{Eq:generalized problem} as \ahmadreza{the following} optimization problem:
\begin{equation}\label{Eq:biconvex reformulation}
    \begin{aligned}
    \max_{\decision\in\mathbb{R}^{\ahmadreza{C}} \atop {\tau\in\mathbb{R}^{I\times N\times \probNum}}}&\enskip \sum_{i\in\mathcal{I}}{\sum_{n\in\mathcal{N}} \sum_{\probInd=1}^\probNum w_{\probInd n}} \decision_i\tau_{i{n\probInd}}\\
    \mbox{s.t.} &\enskip \denuminator_{i{n\probInd}}(\decision)\tau_{i{n\probInd}} \ahmadreza{\leq} 1,& \forall i\in\mathcal{I},\;{ n\in\mathcal{N},\;\probInd=1,...,\probNum}, \\
    &\enskip A\decision\geq b,\\
    &\enskip \tau_{i{n\probInd}},\decision_i\geq 0,&\forall i\in\mathcal{I},\;{ n\in\mathcal{N},\;\probInd=1,...,\probNum}.
    \end{aligned}
\end{equation}
It is clear that \eqref{Eq:generalized problem} and \eqref{Eq:biconvex reformulation} are equivalent as $\denuminator_{i{n\probInd}}(\decision)$ is a positive function \ahmadreza{when $\decision \geq0$, for any $i\in\mathcal{I},$ $n\in\mathcal{N}$, and $\probInd=1,...,\probNum$. Inspired by \cite{zhen2021extension} and the fact that $ \frac{1}{\denuminator_{in\probInd}(\decision)}$ and hence  $\tau_{in\probInd}$ are positive, we can introduce a new variable $\product_{ijn\probInd}=p_j \tau_{in\probInd}$ and rewrite \eqref{Eq:biconvex reformulation} as
\begin{subequations}\label{Eq:biconvex-perspective reformulation}
    \begin{align}
    \max_{\decision\in\mathbb{R}^{\ahmadreza{C}} \atop {\tau\in\mathbb{R}^{I\times N\times \probNum} \atop \product\in\mathbb{R}^{I\times \ahmadreza{\mathcal{C}}\times N\times \probNum}}}&\enskip \sum_{i\in\mathcal{I}}{\sum_{n\in\mathcal{N}} \sum_{\probInd=1}^\probNum w_{\probInd n}} \product_{iin\probInd}\label{eq:obj_perspective}\\
    \mbox{s.t.} &\enskip \denuminator_{i{n\probInd}}(\frac{\product_{i:n\probInd}}{\tau_{i{n\probInd}}})\tau_{i{n\probInd}} \leq 1,& \forall i\in\mathcal{I},\; n\in\mathcal{N},\;\probInd=1,...,\probNum, \label{eq:const_convex}\\
    &\enskip A\decision\geq b,\\
    & \enskip \product_{ijn\probInd}=p_j \tau_{in\probInd}&\forall i\in\mathcal{I},\; j\in\mathcal{C},\;{ n\in\mathcal{N},\;\probInd=1,...,\probNum},\label{eq:const_bilinear}\\
    &\enskip \product_{ijn\probInd},\decision_j\geq 0,&\forall i\in\mathcal{I},\;j\in\mathcal{C},\;{ n\in\mathcal{N},\;\probInd=1,...,\probNum},\\
    &\enskip \tau_{i{n\probInd}}>0&\forall i\in\mathcal{I},\;{ n\in\mathcal{N},\;\probInd=1,...,\probNum},
    \end{align}
\end{subequations}
where $\product_{i:n\probInd}\in\mathbb{R}^{C}$ is the vector containing $\product_{ijn\probInd}$ for $j\in\mathcal{C}$. }

\ahmadreza{Problem \eqref{Eq:biconvex-perspective reformulation} belongs to the class of biconvex optimization problems, as  it contains functions that are convex in $\decision$, and $(\product,\tau)$ but not in $(\decision,\product,\tau)$. To better see this, let us fix $i\in\mathcal{I}, \; j\in\mathcal{C},\;{ n\in\mathcal{N},\;\probInd=1,...,\probNum}$. The function $\denuminator_{{in\probInd}}(\decision)$ is convex as it is the summation of exponential functions with linear exponents in $\decision$. As shown in Section 3.2.6 of \cite{boyd2004convex}, $\denuminator_{i{n\probInd}}(\frac{\product_{i:n\probInd}}{\tau_{i{n\probInd}}})\tau_{i{n\probInd}}$ is the perspective map of the convex function $\denuminator_{{in\probInd}}(\product_{i:n\probInd})$, which is convex in $(\product_{i:n\probInd},\tau_{i{n\probInd}})$. Finally, the term $p_j \tau_{in\probInd}$ appearing in \eqref{eq:const_bilinear} is bilinear. Therefore, \eqref{Eq:biconvex-perspective reformulation} is a biconvex optimization problem. }

\ahmadreza{The main challenge in solving \eqref{Eq:biconvex-perspective reformulation} is on dealing with the bilinear terms in constraint \eqref{eq:const_bilinear}.} 
\ahmadreza{In this section, we use McCormick relaxation \citep{mccormick1976computability} to construct a convex optimization problem that approximates \eqref{Eq:biconvex-perspective reformulation} from above. More specifically, we obtain an upper bound by solving the following optimization problem:}

 \begin{subequations}\label{Eq:Linear approximation}
    \begin{align}
     \max_{\decision\in\mathbb{R}^{\ahmadreza{C}} \atop {{\tau\in\mathbb{R}^{I\times N\times \probNum}} \atop {W\in\mathbb{R}^{I\times I\times N\times \probNum}}}}&\enskip \sum_{i\in\mathcal{I}}{\sum_{n\in\mathcal{N}} \sum_{\probInd=1}^\probNum w_{\probInd n}} W_{ii{n\probInd}}\\
    \mbox{s.t.}&\enskip A\decision\geq b,\;\\
    &\enskip \ahmadreza{\denuminator_{i{n\probInd}}(\frac{\product_{i:n\probInd}}{\tau_{i{n\probInd}}})\tau_{i{n\probInd}} \leq 1},& \forall i\in\mathcal{I},\; n\in\mathcal{N},\;\probInd=1,...,\probNum,\label{eq:const_a convex function}\\
     &\enskip AW_{i:{n\probInd}}\geq b\tau_{i{n\probInd}},& {\forall i\in\mathcal{I},\atop {\forall n\in\mathcal{N},\;\probInd=1,...,\probNum },}\label{Eq: linearization of (Ax-b)tau}\\
    &\enskip LB_{\tau_{i{n\probInd}}}\left( A\decision- b \right) \leq AW_{i:{n\probInd}}- b\tau_{i{n\probInd}},&{\forall i\in\mathcal{I},\atop {\forall n\in\mathcal{N},\;\probInd=1,...,\probNum },}\label{eq:linearization 1 from Zhan,Marandi...}\\
     &\enskip  AW_{i:{n\probInd}}- b\tau_{i{n\probInd}}\leq UB_{\tau_{i{n\probInd}}}\left( A\decision- b \right),&{\forall i\in\mathcal{I},\atop {\forall n\in\mathcal{N},\;\probInd=1,...,\probNum },} \label{eq:linearization 2 from Zhan,Marandi...}\\
    &\enskip W_{ij{n\probInd}}\geq LB_{\tau_{i{n\probInd}}}\decision_j+\tau_{i{n\probInd}} LB_{\decision_j}-LB_{\tau_{i{n\probInd}}}LB_{\decision_j},&{\forall i\in\mathcal{I},\; \ahmadreza{j\in\mathcal{C},}\atop {\forall n\in\mathcal{N},\;\probInd=1,...,\probNum },}\label{Eq:MaCormick 1}\\
     &\enskip W_{ij{n\probInd}}\geq UB_{\tau_{i{n\probInd}}}\decision_j+\tau_{i{n\probInd}} UB_{\decision_j}-UB_{\tau_{i{n\probInd}}}UB_{\decision_j},&{\forall i\in\mathcal{I},\; \ahmadreza{j\in\mathcal{C},}\atop {\forall n\in\mathcal{N},\;\probInd=1,...,\probNum },}\label{Eq:MaCormick2}\\
   &\enskip W_{ij{n\probInd}}\leq UB_{\tau_{i{n\probInd}}}\decision_j+\tau_{i{n\probInd}} LB_{\decision_j}-UB_{\tau_{i{n\probInd}}}LB_{\decision_j},&{\forall i,\in\mathcal{I},\; \ahmadreza{j\in\mathcal{C},}\atop {\forall n\in\mathcal{N},\;\probInd=1,...,\probNum },}\label{Eq:MaCormick 3}\\
     &\enskip      W_{ij{n\probInd}}\geq LB_{\tau_{i{n\probInd}}}\decision_j+\tau_{i{n\probInd}} UB_{\decision_j}-LB_{\tau_{i{n\probInd}}}UB_{\decision_j},&{\forall i,\in\mathcal{I},\; \ahmadreza{j\in\mathcal{C},}\atop {\forall n\in\mathcal{N},\;\probInd=1,...,\probNum },}\label{Eq:MaCormick 4}\\    
     &\enskip LB_\tau\leq \tau \leq UB_\tau,\\
    &\enskip LB_\decision\leq \decision\leq UB_\decision.
    \end{align}
\end{subequations}
where $LB_\decision,UB_\decision\in\mathbb{R}^{\ahmadreza{C}}$ are the vectors containing component-wise lower and upper bounds of $\decision$, \ahmadreza{and}  $LB_\tau,UB_\tau\in\mathbb{R}^{{I\times N\times \probNum }}$ contain the \ahmadreza{positive} component-wise lower and upper bounds of $\tau$, respectively.
Problem \eqref{Eq:Linear approximation} is constructed by \ahmadreza{convexification} of the \ahmadreza{biconvex} optimization problem equivalent to \eqref{Eq:biconvex-perspective reformulation} including some redundant constraints.  
Constraint \eqref{Eq: linearization of (Ax-b)tau} linearizes the redundant constraint $(A\decision-b)\tau_{i{n\probInd}}\geq 0,$ for $i\in\mathcal{I}, { n\in\mathcal{N},\;\probInd=1,...,\probNum }.$  Constraints \eqref{eq:linearization 1 from Zhan,Marandi...} and \eqref{eq:linearization 2 from Zhan,Marandi...} are the constraints proposed by \cite{zhen2018disjoint} to tighten the linear relaxation. Constraints \eqref{Eq:MaCormick 1}, \eqref{Eq:MaCormick2}, \eqref{Eq:MaCormick 3}, and \eqref{Eq:MaCormick 4} are obtained by using McCormick relaxation \citep{mccormick1976computability}. Therefore, the objective value of \eqref{Eq:Linear approximation} is an upper bound on the objective value of \eqref{Eq:biconvex-perspective reformulation} and hence \eqref{Eq:biconvex reformulation}. 

\begin{remark}\label{remark:finding the l and u of tau}
To construct \eqref{Eq:Linear approximation}, we need to compute $LB_\tau$ and $UB_\tau$. Since in the optimal solution $(\tau^*,\decision^*)$ of \eqref{Eq:biconvex reformulation}, we have $\tau^*_{i{n\probInd}}=\frac{1}{f_{i{n\probInd}}(\decision^*)}$, we can compute $\frac{1}{LB_{\tau_{i\ahmadreza{n\probInd}}}}$ by solving  
\begin{equation}\label{eq: lowerbound of tau}
\begin{aligned}
\max_{\decision\in\mathbb{R}^{\ahmadreza{C}}}&\enskip {\denuminator_{i{n\probInd}}(\decision)}\\
 \mbox{s.t.}&\enskip A\decision\geq b, \;\decision \geq 0. 
\end{aligned}
\end{equation}
Since $\denuminator_{i{n\probInd}}(\decision)$ is a convex function, the above optimization problem is known to be $\mathcal{NP}-$hard ({even checking local optimality for a quadratic objective function is $\mathcal{NP}-$hard \citep{pardalos1988checking}}). However, we are not looking for an exact optimal value 
but for an upper bound. \ahmadreza{By definition, we have $\denuminator_{i{\probInd n}}(\decision)={\sum_{j\in \mathcal{C}} e^{V_{jn}(\ahmadreza{b_{\probInd}})-V_{in}(\ahmadreza{b_{\probInd}})}}.$ So, we know
$$
\begin{aligned}
\max_{\decision\in\mathbb{R}^{\ahmadreza{C}}}&\enskip {\denuminator_{i{n\probInd}}(\decision)}\\
 \mbox{s.t.}&\enskip A\decision\geq b, \;\decision \geq 0,
\end{aligned} 
\enskip
\leq 
\enskip
\begin{aligned}
\sum_{j\in \mathcal{C}}\max_{\decision\in\mathbb{R}^{\ahmadreza{C}}}&\enskip e^{V_{jn}(\ahmadreza{b_{\probInd}})-V_{in}(\ahmadreza{b_{\probInd}})}\\
 \mbox{s.t.}&\enskip A\decision\geq b, \;\decision \geq 0,
\end{aligned} 
\enskip 
= 
\enskip 
\sum_{j\in \mathcal{C}}e^ {\left(\begin{aligned}
\max_{\decision\in\mathbb{R}^{\ahmadreza{C}}}&\enskip {V_{jn}(\ahmadreza{b_{\probInd}})-V_{in}(\ahmadreza{b_{\probInd}})}\\
 \mbox{s.t.}&\enskip A\decision\geq b, \;\decision \geq 0,
\end{aligned} \right)} .
$$
So, by solving some linear optimization problems, we can obtain a lower bound on $\tau.$
}
 Also, to compute $\frac{1}{UB_{\tau_{i\ahmadreza{n\probInd}}}}$, we solve
$$
\begin{aligned}
 \min_{\decision\in\mathbb{R}^{\ahmadreza{C}}}&\enskip {\denuminator_{i\ahmadreza{n\probInd}}(\decision)}\\
 \mbox{s.t.}&\enskip A\decision\geq b,\; \decision \geq 0, 
\end{aligned}
$$
which \ahmadreza{is a convex optimization problem; hence} can be solved efficiently. 
\hfill \Halmos
\end{remark}

Hitherto, we have provided a method to obtain a ``good'' feasible solution (Section \ref{Sec:lower bound}) and an optimization problem to provide an upper bound on the objective value of \eqref{Eq:generalized problem} (Section \ref{Sec:upper bound}). {We emphasize that the $p$-part of the solution obtained from \eqref{Eq:Linear approximation} can also be used to find a lower bound on \eqref{Eq:biconvex reformulation}. So, we also check the quality of such solutions.}

In the next section, we \ahmadreza{show} how a branching technique can tighten the gap between the lower and upper bounds. 
\subsection{Branching method}\label{Sec:new branching method}

In this section, we use a typical branching method in continuous optimization. In this method, we first choose the branching variable $i$ and then \ahmadreza{split} its feasible interval into two intervals \citep{misener2014antigone,FLOUDAS20051185,akrotirianakis2004computational}. The index $i$ is usually the dimension where the feasible region has \ahmadreza{its} largest length. Mathematically speaking, 
$$
i \in \argmax\{UB_{\decision_j}-LB_{\decision_j}:\;j\in\ahmadreza{\mathcal{C}}\}.
$$
Then, the partitions are
$$
S_1=\left\{\decision\in\mathbb{R}^{\ahmadreza{C}}:\; LB_{\decision_j}\leq p_i\leq \frac{UB_{\decision_j}+LB_{\decision_j}}{2} \right\}\cap S,\; S_2=\left\{\decision\in\mathbb{R}^{\ahmadreza{C}}:\; \frac{UB_{\decision_j}+LB_{\decision_j}}{2}\leq p_i\leq LB_{\decision_j} \right \}\cap S.
$$

As one can notice, such branching methods result in binary trees, as in each iteration\ahmadreza{,} we only have two branches. Furthermore, in each iteration\ahmadreza{,} the volume of the feasible region is getting smaller. Therefore, we can have the following convergence result.

\begin{theorem}\label{TH:Assysmtotically convergence}
Let us denote by $S$ the feasible region of \eqref{Eq:generalized problem}, and its partitions $S^m_k$, $k=1,...,K^m$, in the $m^{\mbox{th}}$ iteration. Also, let us denote by $\mathcal{B}_\ahmadreza{t}(\decision)$ a hyperball with the center $\decision$ and radius $\ahmadreza{t}$.  Let $opt^m$ be the upper bound obtained in the $m^{\mbox{th}}$ iteration. Set 
$$
\ahmadreza{t}^m:=\max_{k=1,...,K^m}\left\{\min\left\{\ahmadreza{t}: S^m_k\subseteq\mathcal{B}_\ahmadreza{t}(\decision), \mbox{for some }\decision\in S^m_k\right\}\right\}.
$$
In other words, $\ahmadreza{t}^m$ is the maximum radius of the smallest hyperball among those covering the partitions $S^m_k$. If $\ahmadreza{t}_m\rightarrow 0$ as $m$ tends to $+\infty$, then  $opt^m\searrow opt$, meaning the sequence of upper bounds asymptotically converges to the optimal value of \eqref{Eq:generalized problem}.
\end{theorem}

\begin{proof}{Proof.}
\ahmadreza{
   Without loss of generality, we assume that $b_\probInd$ can be written as $[b^p_\probInd \; b^q_\probInd]$, where $b^p_\probInd, b^q_\probInd \in \mathbb{R}^{C\times N}$, for any $\probInd=1,...,\probNum$. Let us set $M_{n\probInd} = \max_{i,j\in \mathcal{C}}\left\lVert \left[b^p_{jn\probInd} \atop -b^p_{in\probInd} \right] \right\rVert_2$, for any $n\in\mathcal{N}$ and $\probInd = 1,...,\probNum$. Let us assume that $\decision \in \mathcal{B}_t(\decision^0)$, for a given feasible $\decision^0$ and $t>0$. We see from Remark \ref{remark:finding the l and u of tau} that 
    $$
   \sum_{j\in \mathcal{C}}e^ {
\min_{\decision \in \mathcal{B}_\ahmadreza{t}(\decision^0)} \{V_{jn}(\ahmadreza{b_{\probInd}})-V_{in}(\ahmadreza{b_{\probInd}})\}}  \leq \frac{1}{\tau_{in\probInd}} \leq \sum_{j\in \mathcal{C}}e^ {
\max_{\decision \in \mathcal{B}_\ahmadreza{t}(\decision^0)} \{V_{jn}(\ahmadreza{b_{\probInd}})-V_{in}(\ahmadreza{b_{\probInd}})\}}.
    $$
    For each $j\in\mathcal{C}$, we know 
    \begin{eqnarray}
    \max_{\decision \in \mathcal{B}_\ahmadreza{t}(\decision^0)} \{V_{jn}(\ahmadreza{b_{\probInd}})-V_{in}(\ahmadreza{b_{\probInd}})\} &=&  q_{in}(b^q_\probInd) -q_{jn}(b^q_\probInd) + \max_{\decision \in \mathcal{B}_\ahmadreza{t}(\decision^0)} \left[b^p_{jn\probInd} \atop -b^p_{in\probInd} \right]^T\left[p_j \atop p_i \right]\nonumber \\ 
    &=&  q_{in}(b^q_\probInd) -q_{jn}(b^q_\probInd) + \left[b^p_{jn\probInd} \atop -b^p_{in\probInd} \right]^T\left[\decision^0_j \atop \decision^0_i \right] +t\left\lVert  \left[b^p_{jn\probInd} \atop -b^p_{in\probInd} \right]\right\rVert_2 \nonumber\\
    &\leq&  q_{in}(b^q_\probInd) -q_{jn}(b^q_\probInd) + \left[b^p_{jn\probInd} \atop -b^p_{in\probInd} \right]^T\left[\decision^0_j \atop \decision^0_i \right] +tM_{n\probInd},
    \end{eqnarray}
    and, analogously,
    $$
    \min_{\decision \in \mathcal{B}_\ahmadreza{t}(\decision^0)} \{V_{jn}(\ahmadreza{b_{\probInd}})-V_{in}(\ahmadreza{b_{\probInd}})\} \geq q_{in}(b^q_\probInd) -q_{jn}(b^q_\probInd) + \left[b^p_{jn\probInd} \atop -b^p_{in\probInd} \right]^T\left[\decision^0_j \atop \decision^0_i \right] -tM_{n\probInd}.
    $$
    So, we have 
    $$
   \frac{e^{-t M_{n\probInd}}}{f_{in\probInd}(\decision^0)}  \leq \tau_{in\probInd} \leq \frac{e^{t M_{n\probInd}}}{f_{in\probInd}(\decision^0)}.
    $$
    Therefore, for given $i\in\mathcal{I}$, $n\in\mathcal{N}$, and $\probInd=1,...,\probNum$,  when $t\rightarrow 0$, the differences of the lower and upper bounds of $\decision_{i}$ and $\tau_{in\probInd}$, tends to zero. Furthermore, $opt^m$ is a non-increasing sequence. Hence, if $t_m\rightarrow 0$ as $m$ tends to $+\infty$, then $opt^m \searrow opt$ \citep{mccormick1976computability}.
    }
    \hfill \Halmos
\end{proof}
\ahmadreza{The main assumption in Theorem \ref{TH:Assysmtotically convergence} is that the radius of the minimum ball covering the partitions goes to $0$. This assumption is satisfied by the partitioning method used in the algorithm. Therefore, the theorem guarantees that the algorithm achieves an $\epsilon-$approximation of the optimal value of \eqref{Eq:ourproblem} within a finite number of steps for any given $\epsilon >0$.
 }
With this assumption, we can obtain multiple optimal solutions from the trust-region approach as well as the $p-$part of the solution of \eqref{Eq:Linear approximation}. So, the proposed algorithm is able to provide the user a collection of optimal solutions, which is important in pricing problems. 

In the next section, we show how the algorithm efficiently works on several instances.

\section{Numerical experiments}
In this section, we discuss the effectiveness of our method \ahmadreza{in solving} static pricing problems under a discrete mixed logit model. We refer to our method as \ahmadreza{CoBiT}, as it is based on \ahmadreza{\textbf{Co}vexification} of a \textbf{bi}convex optimization and \textbf{t}rust-region algorithm. To have a better understanding of how \ahmadreza{CoBiT} works, we start by presenting its logic on \ahmadreza{an} illustrative example in Section \ref{Sec:numerical_illustration}.  We then analyze \ahmadreza{CoBiT} from an algorithmic perspective in Section \ref{Sec:lidata} and Section \ref{MSdata} by comparing its performance to state-of-the-art nonlinear optimizers and the two closest algorithms available in the literature: the algorithm 2 of \cite{li2019product} and the algorithm of \cite{GeerBranch}. Finally, Section \ref{Sec:parking choice model} investigates \ahmadreza{CoBiT} from a problem definition perspective by comparing the results obtained with pricing problems that cannot handle price sensitivity parameters
that differ among products and customer segments. To do so, we use an additional case study, inspired by a real parking choice model. \color{black}

The numerical results of this work were carried out on a Laptop featuring 4 processors   2.60 GHz and 8.00 GB RAM running Julia \ahmadreza{1.8.0} \citep{Julia}. We use JuMP \ahmadreza{1.7.0} \citep{JuMP} to pass Linear Optimization problems to Gurobi 9.5.1. 
\ahmadreza{We solve \eqref{Eq:Linear approximation} by reformulating it in a conic form and using Mosek Optimization Tool 10.0.}
We compare the performance of our method with SCIP 7.0 \citep{gamrath2020scip}, Couenne \citep{belotti2009couenne}, BARON 21.1.13 \citep{sahinidis1996baron}, and ANTIGONE 37.1 \citep{misener2014antigone}. We have passed the problems to SCIP, BARON, and ANTIGONE using GAMS modeling language version 37.1.0. For Couenne, we use JuMP to pass the problem to the solver.
 

\subsection{Illustration of \ahmadreza{CoBiT}}\label{Sec:numerical_illustration}
To solve the optimization problem, \ahmadreza{CoBiT} starts by partitioning the feasible region into two parts. The trust-region algorithm, Algorithm \ref{alg:trust region}, is employed to obtain new local solutions in each of the rectangles (Figure\ref{fig:first iteration}).  

Then, \ahmadreza{CoBiT} further partitions each of the rectangles (Figure \ref{fig:second iteration}). After using the trust region in the top left square, \ahmadreza{CoBiT} finds the same local solution as in the previous iteration. In the top right square, \ahmadreza{CoBiT} finds a new solution (red square in the figure). For both two bottom squares, \ahmadreza{CoBiT} detects that no effort is needed there as the upper bounds on these regions are, at most, as high as the objective value of the best-obtained solution. So, in the next iteration, the top squares are further partitioned, and local solutions are obtained (blue dots in Figure\ref{fig:third iteration}).


\ahmadreza{CoBiT} continues the procedure until either the time limit is reached or the objective value of the best obtained feasible solution does not deviate by more than $10^{-5}$ from the upper bound obtained by the linearization.  
\begin{figure}[h!]
\centering
    \begin{subfigure}{0.3\linewidth}
        \centering
\begin{tikzpicture}[scale=0.8]
\filldraw[black, very thick] (3,2.13) ellipse (3pt and 1.5pt);
\filldraw[black, very thick] (1.5,1.9) ellipse (3pt and 1.5pt);
\draw[black,thick] (1.5,0)--(1.5,3);
\filldraw[fill=blue!40!white, draw=black, opacity=0.2] (0,0) rectangle (3,3);
\end{tikzpicture} 
\subcaption{First iteration}
    \label{fig:first iteration}
    \end{subfigure}
     \begin{subfigure}{0.3\linewidth}
        \centering
\begin{tikzpicture}[scale=0.8]
\filldraw[black, very thick] (3,2.13) ellipse (3pt and 1.5pt);
\filldraw[black, very thick] (1.5,1.9) ellipse (3pt and 1.5pt);
\draw[black,thick] (1.5,0)--(1.5,3);
\draw[black,thick] (0,1.5)--(3,1.5);
\filldraw[red, very thick] (2.1,2.13) +(-2pt,-2pt) rectangle +(2pt,2pt);
%
\filldraw[fill=blue!40!white, draw=black, opacity=0.2] (0,0) rectangle (3,3);
\end{tikzpicture} 
\subcaption{Second iteration}
    \label{fig:second iteration}
    \end{subfigure}
        \begin{subfigure}{0.3\linewidth}
        \centering
\begin{tikzpicture}[scale=0.8]
\filldraw[black, very thick] (3,2.13) ellipse (3pt and 1.5pt);
\filldraw[black, very thick] (1.5,1.9) ellipse (3pt and 1.5pt);
\draw[black,thick] (1.5,0)--(1.5,3);
\draw[black,thick] (0,1.5)--(3,1.5);
\filldraw[red, very thick] (2.1,2.13) +(-2pt,-2pt) rectangle +(2pt,2pt);
\filldraw[blue, very thick] (2.25,2.14) circle (0.1);
\filldraw[blue, very thick] (0.75,1.5) circle (0.1);

\draw[black,thick] (2.25,1.5)--(2.25,3.0);
\draw[black,thick] (0.75,1.5)--(0.75,3.0);
\filldraw[fill=blue!40!white, draw=black, opacity=0.2] (0,0) rectangle (3,3);
\end{tikzpicture} 
\subcaption{Third iteration}
    \label{fig:third iteration}
    \end{subfigure}
\caption{An illustration of first three iterations of \ahmadreza{CoBiT}}\color{black}\label{Fig: illustration of the numerical experimentsN=10}
\end{figure}

\subsection{\ahmadreza{CoBiT} versus \cite{li2019product}} \label{Sec:lidata}
In this section, we apply \ahmadreza{CoBiT} to the Intel Corporation case presented in \cite{li2019product}.

\subsubsection{The Intel Corporation case} 
The authors assume that the Intel market is divided into a finite number of segments and \ahmadreza{that a multinomial logit model characterizes the product demand in each segment.} Unlike us, the authors do not allow any \ahmadreza{customer} specific variables in the utility functions, and their pricing optimization problem is therefore given by:
 
\begin{equation}\label{Eq:MixedMNL}
\begin{aligned}
\max_{p\in\mathbb{R}^I}& \enskip \sum_{i\in\mathcal{I}} p_i P_{i}, & \\
\mbox{s.t.}&\enskip P_{i} = \sum_{k=1}^{K} w_{k} \frac{e^{q_{ik}+\beta_{ik} p_i}}{\sum_{j\in \mathcal{I}} e^{q_{jk}+\beta_{jk} p_j}}, &\forall i \in \mathcal{I}.
\end{aligned}
\end{equation}
\vspace{0.15cm}

The authors apply this model to Intel’s microprocessor stock-keeping units (SKUs) used in computer servers.
Intel customers are categorized into $K=7$ segments, corresponding to the seven groupings used by Intel’s sales division. The weights $\boldsymbol{w_{k}}$, computed based on historical purchasing volumes, are provided in Table \ref{tab:inteldata3}.  

\begin{table}[h]
    \centering
    \begin{tabular}{|c|c|c|c|c|c|c|c|}
     \hline
      k  &  1 & 2 & 3 & 4 & 5 & 6 & 7 \\
       \hline
      $w_k$ & 0.0753 & 0.1126 & 0.1285 & 0.1180 & 0.0859 & 0.2842 & 0.1953 \\
      \hline
    \end{tabular}
    \caption{$\boldsymbol{w_{k}}$  values from  \cite{li2019product}}
    \label{tab:inteldata3}
\end{table}

Sales data of the first three generations of products (13 SKUs) were used to parameterize the demand model, while the fourth generation of products, representing 3 SKUs ($\mathcal{I}=\{1,2,3\}$), were used to test the demand model. The segment-specific coefficients are provided in Tables \ref{tab:inteldata1} and \ref{tab:inteldata2}. We refer to \cite{li2019product}, the original work, for more details on data fitting and parameterization. Prices were then optimized for the three SKUs of this fourth generation of products.

\begin{table}[h]
    \centering
    \begin{tabular}{|c|c|c|c|c|c|c|c|}
     \hline
       \multicolumn{8}{|c|}{$k$}  \\
       \hline
      $i$   &  1 & 2 & 3 & 4 & 5 & 6 & 7 \\
      \hline
      1 & -1.0334 & 3.2480 & -0.9336 & 1.7094 & 0.4187 & -0.8904 & -0.9804\\
      \hline
      2 & 0.7840  & 4.7161 & -0.3438 & 1.8777 & 2.1771 & -0.4310 & -0.4907 \\
      \hline
      3 & 6.0054 & 3.8771 & 1.3506 & 2.3611 & 1.1723 & 0.8889 & 0.9163 \\ 
      \hline
    \end{tabular}
    \caption{$\boldsymbol{q_{ik}}$  values from  \cite{li2019product}}
    \label{tab:inteldata1}
\end{table}

\begin{table}[h]
    \centering
    \begin{tabular}{|c|c|c|c|c|c|c|c|}
     \hline
       \multicolumn{8}{|c|}{$k$}  \\
       \hline
      $i$   &  1 & 2 & 3 & 4 & 5 & 6 & 7 \\
      \hline
      1 & -0.00416 &  -0.01840 &  -0.00525 &  -0.01165 &  -0.01015 &  -0.00325 &  -0.00331\\
      \hline
      2 &  -0.00312 & -0.01354 & -0.00394 & -0.00874 & -0.00639 & -0.00244 & -0.00248\\
      \hline
      3 & -0.00181 & -0.00744 & -0.00229 & -0.00508 & -0.00167 & -0.00142 & -0.00144 \\ 
      \hline
    \end{tabular}
    \caption{$\boldsymbol{\beta_{ik}}$  values from  \cite{li2019product}}
    \label{tab:inteldata2}
\end{table}

\subsubsection{Numerical results}
Using the provided data, Algorithm 2 in \cite{li2019product} finds the local solution $[608.2695 \;\;365.079 \;\;1,209.09]$ within 0.079 seconds with the objective value of $362.3389$. \ahmadreza{CoBiT} can find this solution \ahmadreza{and prove optimality in 105.78} seconds. As explained in \cite{li2019product}, their Algorithm 2 performs very well ``when the degree of segment asymmetry is sufficiently small, because then the total profit is concave'', which is the case for this instance. 

To see the effect of this algorithm in other cases and compare it fairly with \ahmadreza{CoBiT}, we generate 10 instances randomly. More specifically, we choose $q_{ik}$ from a uniform distribution over $[-5,5]$ and $\beta_{ik}$ from a uniform distribution over $[-5,-0.025]$. Table \ref{Table:Li et al vs libit} provides the result of applying \ahmadreza{CoBiT} and Algorithm 2 in \cite{li2019product} on these 10 instances.

\begin{table}[h]
\centering
\begin{tabular}{ccccccc}
\hline
\multicolumn{1}{|c|}{}                                    & \multicolumn{1}{c|}{}                         & \multicolumn{1}{c|}{}                                   & \multicolumn{1}{c|}{}                                & \multicolumn{3}{c|}{Price}                                                           \\ \cline{5-7} 
\multicolumn{1}{|c|}{\multirow{-2}{*}{{Instance}}} & \multicolumn{1}{c|}{\multirow{-2}{*}{Method}} & \multicolumn{1}{c|}{\multirow{-2}{*}{Expected revenue}} & \multicolumn{1}{c|}{\multirow{-2}{*}{Solution time}} & \multicolumn{1}{c|}{$p_1$}  & \multicolumn{1}{c|}{$p_2$}  & \multicolumn{1}{c|}{$p_3$}  \\ \hline
\multicolumn{1}{|c|}{}                                                            & \multicolumn{1}{c|}{\ahmadreza{CoBiT}}                    & \multicolumn{1}{c|}{ {1.40}}                             & \multicolumn{1}{c|}{\ahmadreza{672.88}}                                & \multicolumn{1}{c|}{ 4.\ahmadreza{19}}  & \multicolumn{1}{c|}{ \ahmadreza{7.81}}  & \multicolumn{1}{c|}{ 3.\ahmadreza{13}}  \\ 
\multicolumn{1}{|c|}{\multirow{-2}{*}{1}}                                         & \multicolumn{1}{c|}{Algorithm 2 \citep{li2019product}}                 & \multicolumn{1}{c|}{1.40}                             & \multicolumn{1}{c|}{2.32}                            & \multicolumn{1}{c|}{4.23}  & \multicolumn{1}{c|}{8.04}  & \multicolumn{1}{c|}{2.92}  \\ \hline
\multicolumn{1}{|c|}{}                                                            & \multicolumn{1}{c|}{{\ahmadreza{CoBiT}}}           & \multicolumn{1}{c|}{\textbf{1.41}}                    & \multicolumn{1}{c|}{\ahmadreza{71.67}}                       & \multicolumn{1}{c|}{5.7}  & \multicolumn{1}{c|}{13.50} & \multicolumn{1}{c|}{10.50} \\ 
\multicolumn{1}{|c|}{\multirow{-2}{*}{{2}}}                                & \multicolumn{1}{c|}{Algorithm 2 \citep{li2019product}}                 & \multicolumn{1}{c|}{1.37}                             & \multicolumn{1}{c|}{5.77}                            & \multicolumn{1}{c|}{5.71}  & \multicolumn{1}{c|}{10.46} & \multicolumn{1}{c|}{4.61}  \\ \hline
\multicolumn{1}{|c|}{}                                                            & \multicolumn{1}{c|}{\ahmadreza{CoBiT}}                    & \multicolumn{1}{c|}{\textbf{4.46}}                              & \multicolumn{1}{c|}{\ahmadreza{3.94}}                                & \multicolumn{1}{c|}{68.\ahmadreza{60}} & \multicolumn{1}{c|}{2.\ahmadreza{50}}  & \multicolumn{1}{c|}{17.\ahmadreza{40}} \\ 
\multicolumn{1}{|c|}{\multirow{-2}{*}{3}}                                         & \multicolumn{1}{c|}{Algorithm 2 \citep{li2019product}}                 & \multicolumn{1}{c|}{0.554}                              & \multicolumn{1}{c|}{3.11}                            & \multicolumn{1}{c|}{1.0}     & \multicolumn{1}{c|}{1.36}  & \multicolumn{1}{c|}{0.85}  \\ \hline
\multicolumn{1}{|c|}{}                                                            & \multicolumn{1}{c|}{\ahmadreza{CoBiT}}                    & \multicolumn{1}{c|}{0.65}                             & \multicolumn{1}{c|}{\ahmadreza{313.36}}                                & \multicolumn{1}{c|}{1.38}  & \multicolumn{1}{c|}{3.\ahmadreza{00}}  & \multicolumn{1}{c|}{1.39}  \\ 
\multicolumn{1}{|c|}{\multirow{-2}{*}{4}}                                         & \multicolumn{1}{c|}{Algorithm 2 \citep{li2019product}}                 & \multicolumn{1}{c|}{0.65}                             & \multicolumn{1}{c|}{1.95}                            & \multicolumn{1}{c|}{1.38}  & \multicolumn{1}{c|}{3.11}  & \multicolumn{1}{c|}{1.39}  \\ \hline
\multicolumn{1}{|c|}{}                                                            & \multicolumn{1}{c|}{\ahmadreza{CoBiT}}                    & \multicolumn{1}{c|}{0.80}                                & \multicolumn{1}{c|}{\ahmadreza{267.89}}                                & \multicolumn{1}{c|}{1.2\ahmadreza{5}}  & \multicolumn{1}{c|}{3.6\ahmadreza{3}}  & \multicolumn{1}{c|}{1.1\ahmadreza{8}}  \\ 
\multicolumn{1}{|c|}{\multirow{-2}{*}{5}}                                         & \multicolumn{1}{c|}{Algorithm 2 \citep{li2019product}}                 & \multicolumn{1}{c|}{0.80}                                & \multicolumn{1}{c|}{2.57}                            & \multicolumn{1}{c|}{0.32}  & \multicolumn{1}{c|}{0.91}  & \multicolumn{1}{c|}{0.29}  \\ \hline
\multicolumn{1}{|c|}{}                                                            & \multicolumn{1}{c|}{\ahmadreza{CoBiT}}                    & \multicolumn{1}{c|}{1.24}                               & \multicolumn{1}{c|}{\ahmadreza{91.17}}                                & \multicolumn{1}{c|}{2.19}  & \multicolumn{1}{c|}{6.1\ahmadreza{3}}  & \multicolumn{1}{c|}{2.04}  \\ 
\multicolumn{1}{|c|}{\multirow{-2}{*}{6}}                                         & \multicolumn{1}{c|}{Algorithm 2 \citep{li2019product}}                 & \multicolumn{1}{c|}{1.24}                               & \multicolumn{1}{c|}{1.59}                            & \multicolumn{1}{c|}{2.19}  & \multicolumn{1}{c|}{6.16}  & \multicolumn{1}{c|}{2.04}  \\ \hline
\multicolumn{1}{|c|}{}                                                            & \multicolumn{1}{c|}{\ahmadreza{CoBiT}}                    & \multicolumn{1}{c|}{2.20}                                & \multicolumn{1}{c|}{\ahmadreza{74.92}}                                & \multicolumn{1}{c|}{2.2\ahmadreza{5}}  & \multicolumn{1}{c|}{0.92}  & \multicolumn{1}{c|}{15.0\ahmadreza{8}} \\ 
\multicolumn{1}{|c|}{\multirow{-2}{*}{7}}                                         & \multicolumn{1}{c|}{Algorithm 2 \citep{li2019product}}                 & \multicolumn{1}{c|}{2.20}                                & \multicolumn{1}{c|}{1.75}                            & \multicolumn{1}{c|}{2.26}  & \multicolumn{1}{c|}{0.92}  & \multicolumn{1}{c|}{15.07} \\ \hline
\multicolumn{1}{|c|}{}                                                            & \multicolumn{1}{c|}{\ahmadreza{CoBiT}}                    & \multicolumn{1}{c|}{0.69}                             & \multicolumn{1}{c|}{\ahmadreza{153.51}}                                & \multicolumn{1}{c|}{0.8\ahmadreza{6}}  & \multicolumn{1}{c|}{\ahmadreza{3.98}}  & \multicolumn{1}{c|}{0.94}  \\ 
\multicolumn{1}{|c|}{\multirow{-2}{*}{8}}                                         & \multicolumn{1}{c|}{Algorithm 2 \citep{li2019product}}                 & \multicolumn{1}{c|}{0.69}                             & \multicolumn{1}{c|}{1.64}                            & \multicolumn{1}{c|}{0.84}  & \multicolumn{1}{c|}{4.33}  & \multicolumn{1}{c|}{0.94}  \\ \hline
\multicolumn{1}{|c|}{}                                                            & \multicolumn{1}{c|}{\ahmadreza{CoBiT}}                    & \multicolumn{1}{c|}{3.66}                              & \multicolumn{1}{c|}{\ahmadreza{19.99}}                                & \multicolumn{1}{c|}{65.\ahmadreza{73}} & \multicolumn{1}{c|}{1.1\ahmadreza{4}}  & \multicolumn{1}{c|}{16.\ahmadreza{32}} \\ 
\multicolumn{1}{|c|}{\multirow{-2}{*}{9}}                                         & \multicolumn{1}{c|}{Algorithm 2 \citep{li2019product}}                 & \multicolumn{1}{c|}{3.66}                              & \multicolumn{1}{c|}{1.77}                            & \multicolumn{1}{c|}{65.57} & \multicolumn{1}{c|}{1.13}  & \multicolumn{1}{c|}{16.29} \\ \hline
\multicolumn{1}{|c|}{}                                                            & \multicolumn{1}{c|}{\ahmadreza{CoBiT}}                    & \multicolumn{1}{c|}{\textbf{0.91}}                             & \multicolumn{1}{c|}{\ahmadreza{40.19}}                                & \multicolumn{1}{c|}{\ahmadreza{9.00}}  & \multicolumn{1}{c|}{13.\ahmadreza{50}}  & \multicolumn{1}{c|}{7.\ahmadreza{20}}  \\ 
\multicolumn{1}{|c|}{\multirow{-2}{*}{10}}                                        & \multicolumn{1}{c|}{Algorithm 2 \citep{li2019product}}                 & \multicolumn{1}{c|}{0.72}                             & \multicolumn{1}{c|}{2.12}                            & \multicolumn{1}{c|}{1.29}  & \multicolumn{1}{c|}{2.16}  & \multicolumn{1}{c|}{1.54}  \\ \hline
\end{tabular}
\caption{Comparison of the solution obtained by the local algorithm in Algorithm 2 \citep{li2019product} and global solver \ahmadreza{CoBiT}. The time inside parentheses is the time took for \ahmadreza{CoBiT} to find an optimal solution.  }\label{Table:Li et al vs libit}
\end{table}

As one can see, for 7 out of 10 randomly generated instances, Algorithm 2 in \cite{li2019product} finds an optimal solution of the problem. However, in the other 3 instances, the optimality gap can go up to $80\%$. 
From the comparison of the optimal prices in these 3 instances, we see the price vector reported by Algorithm 2 is far from the optimal prices. This shows that the local solutions may not necessarily be close to the optimal prices. 


Furthermore, we can see that in \ahmadreza{5} instances (\ahmadreza{2}, 3, 7, 9, and \ahmadreza{10}), the solution time of \ahmadreza{CoBiT} is much smaller than the others. The main reason behind it is the shape of the expected revenue function. If the function is rather flat, then the upper bounds obtained by solving Problem (\ref{Eq:Linear approximation}) may not help much in deciding not to branch some nodes. In other words, when the expected revenue function is flat, then we need to have partitions that are really small to have a good \ahmadreza{convex} approximation. However, in the \ahmadreza{5} instances where \ahmadreza{CoBiT} is fast, the expected revenue function is not very flat, and therefore, the number of nodes that are explored in each iteration of \ahmadreza{CoBiT} \ahmadreza{is} reasonable, which results in a more reasonable solution time. 

\subsection{\ahmadreza{CoBiT} versus \cite{GeerBranch}} 
\label{MSdata}
\ahmadreza{
\begin{table}[h!]
\ahmadreza{
\begin{tabular}{lllllllllllllllllllllllllllllllllllllllllllllllllllllll}
\cline{1-11}
\multicolumn{1}{|l|}{Instance} & \multicolumn{2}{c|}{330}                                   & \multicolumn{2}{c|}{331}                                   & \multicolumn{2}{c|}{332}                                   & \multicolumn{2}{c|}{430}                                  & \multicolumn{2}{c|}{431}                                    &  &  &  &  &  &  &  &  &  &  &  &  &  &  &  &  &  &  &  &  &  &  &  &  &  &  &  &  &  &  &  &  &  &  &  &  &  &  &  &  &  &  &  &  \\ \cline{1-11}
\multicolumn{1}{|l|}{Method}     & \multicolumn{1}{l|}{CoBiT}  & \multicolumn{1}{l|}{GB}      & \multicolumn{1}{l|}{CoBiT}  & \multicolumn{1}{l|}{GB}      & \multicolumn{1}{l|}{CoBiT}  & \multicolumn{1}{l|}{GB}      & \multicolumn{1}{l|}{CoBiT}  & \multicolumn{1}{l|}{GB}     & \multicolumn{1}{l|}{CoBiT}   & \multicolumn{1}{l|}{GB}      &  &  &  &  &  &  &  &  &  &  &  &  &  &  &  &  &  &  &  &  &  &  &  &  &  &  &  &  &  &  &  &  &  &  &  &  &  &  &  &  &  &  &  &  \\ \cline{1-11}
\multicolumn{1}{|l|}{Time (S)}     & \multicolumn{1}{l|}{32.78}  & \multicolumn{1}{l|}{7200}    & \multicolumn{1}{l|}{9.49}   & \multicolumn{1}{l|}{613.57}  & \multicolumn{1}{l|}{17.05}  & \multicolumn{1}{l|}{249.68}  & \multicolumn{1}{l|}{35.59}  & \multicolumn{1}{l|}{62.25}  & \multicolumn{1}{l|}{16.33}   & \multicolumn{1}{l|}{84.34}   &  &  &  &  &  &  &  &  &  &  &  &  &  &  &  &  &  &  &  &  &  &  &  &  &  &  &  &  &  &  &  &  &  &  &  &  &  &  &  &  &  &  &  &  \\ \cline{1-11}
\multicolumn{1}{|l|}{Gap (\%)}    & \multicolumn{1}{l|}{0}      & \multicolumn{1}{l|}{1.37}    & \multicolumn{1}{l|}{0}      & \multicolumn{1}{l|}{0}       & \multicolumn{1}{l|}{0}      & \multicolumn{1}{l|}{0}       & \multicolumn{1}{l|}{0}      & \multicolumn{1}{l|}{0}      & \multicolumn{1}{l|}{0}       & \multicolumn{1}{l|}{-0.03}   &  &  &  &  &  &  &  &  &  &  &  &  &  &  &  &  &  &  &  &  &  &  &  &  &  &  &  &  &  &  &  &  &  &  &  &  &  &  &  &  &  &  &  &  \\ \cline{1-11}
                               &                             &                              &                             &                              &                             &                              &                             &                             &                              &                              &  &  &  &  &  &  &  &  &  &  &  &  &  &  &  &  &  &  &  &  &  &  &  &  &  &  &  &  &  &  &  &  &  &  &  &  &  &  &  &  &  &  &  &  \\ \cline{1-11}
\multicolumn{1}{|l|}{Instance} & \multicolumn{2}{c|}{432}                                   & \multicolumn{2}{c|}{530}                                   & \multicolumn{2}{c|}{531}                                   & \multicolumn{2}{c|}{532}                                  & \multicolumn{2}{c|}{340}                                    &  &  &  &  &  &  &  &  &  &  &  &  &  &  &  &  &  &  &  &  &  &  &  &  &  &  &  &  &  &  &  &  &  &  &  &  &  &  &  &  &  &  &  &  \\ \cline{1-11}
\multicolumn{1}{|l|}{Method}     & \multicolumn{1}{l|}{CoBiT}  & \multicolumn{1}{l|}{GB}      & \multicolumn{1}{l|}{CoBiT}  & \multicolumn{1}{l|}{GB}      & \multicolumn{1}{l|}{CoBiT}  & \multicolumn{1}{l|}{GB}      & \multicolumn{1}{l|}{CoBiT}  & \multicolumn{1}{l|}{GB}     & \multicolumn{1}{l|}{CoBiT}   & \multicolumn{1}{l|}{GB}      &  &  &  &  &  &  &  &  &  &  &  &  &  &  &  &  &  &  &  &  &  &  &  &  &  &  &  &  &  &  &  &  &  &  &  &  &  &  &  &  &  &  &  &  \\ \cline{1-11}
\multicolumn{1}{|l|}{Time (S)}     & \multicolumn{1}{l|}{74.08}  & \multicolumn{1}{l|}{89.49}   & \multicolumn{1}{l|}{268.90} & \multicolumn{1}{l|}{470.47}  & \multicolumn{1}{l|}{526.73} & \multicolumn{1}{l|}{520.81}  & \multicolumn{1}{l|}{107.00} & \multicolumn{1}{l|}{183.39} & \multicolumn{1}{l|}{5.47}     & \multicolumn{1}{l|}{7200}    &  &  &  &  &  &  &  &  &  &  &  &  &  &  &  &  &  &  &  &  &  &  &  &  &  &  &  &  &  &  &  &  &  &  &  &  &  &  &  &  &  &  &  &  \\ \cline{1-11}
\multicolumn{1}{|l|}{Gap (\%)}    & \multicolumn{1}{l|}{0}      & \multicolumn{1}{l|}{0}       & \multicolumn{1}{l|}{0}      & \multicolumn{1}{l|}{0}       & \multicolumn{1}{l|}{0}      & \multicolumn{1}{l|}{0}       & \multicolumn{1}{l|}{0}      & \multicolumn{1}{l|}{-0.01}  & \multicolumn{1}{l|}{0}       & \multicolumn{1}{l|}{0.04}    &  &  &  &  &  &  &  &  &  &  &  &  &  &  &  &  &  &  &  &  &  &  &  &  &  &  &  &  &  &  &  &  &  &  &  &  &  &  &  &  &  &  &  &  \\ \cline{1-11}
                               &                             &                              &                             &                              &                             &                              &                             &                             &                              &                              &  &  &  &  &  &  &  &  &  &  &  &  &  &  &  &  &  &  &  &  &  &  &  &  &  &  &  &  &  &  &  &  &  &  &  &  &  &  &  &  &  &  &  &  \\ \cline{1-11}
\multicolumn{1}{|l|}{Instance} & \multicolumn{2}{c|}{341}                                   & \multicolumn{2}{c|}{342}                                   & \multicolumn{2}{c|}{440}                                   & \multicolumn{2}{c|}{441}                                  & \multicolumn{2}{c|}{442}                                    &  &  &  &  &  &  &  &  &  &  &  &  &  &  &  &  &  &  &  &  &  &  &  &  &  &  &  &  &  &  &  &  &  &  &  &  &  &  &  &  &  &  &  &  \\ \cline{1-11}
\multicolumn{1}{|l|}{Method}     & \multicolumn{1}{l|}{CoBiT}  & \multicolumn{1}{l|}{GB}      & \multicolumn{1}{l|}{CoBiT}  & \multicolumn{1}{l|}{GB}      & \multicolumn{1}{l|}{CoBiT}  & \multicolumn{1}{l|}{GB}      & \multicolumn{1}{l|}{CoBiT}  & \multicolumn{1}{l|}{GB}     & \multicolumn{1}{l|}{CoBiT}   & \multicolumn{1}{l|}{GB}      &  &  &  &  &  &  &  &  &  &  &  &  &  &  &  &  &  &  &  &  &  &  &  &  &  &  &  &  &  &  &  &  &  &  &  &  &  &  &  &  &  &  &  &  \\ \cline{1-11}
\multicolumn{1}{|l|}{Time (S)}     & \multicolumn{1}{l|}{35.6}  & \multicolumn{1}{l|}{169.34}  & \multicolumn{1}{l|}{5.06}    & \multicolumn{1}{l|}{7200}    & \multicolumn{1}{l|}{13.68}  & \multicolumn{1}{l|}{6876.57} & \multicolumn{1}{l|}{4.23}   & \multicolumn{1}{l|}{217.15} & \multicolumn{1}{l|}{38.29}   & \multicolumn{1}{l|}{4717.28} &  &  &  &  &  &  &  &  &  &  &  &  &  &  &  &  &  &  &  &  &  &  &  &  &  &  &  &  &  &  &  &  &  &  &  &  &  &  &  &  &  &  &  &  \\ \cline{1-11}
\multicolumn{1}{|l|}{Gap (\%)}    & \multicolumn{1}{l|}{0}      & \multicolumn{1}{l|}{-0.04}   & \multicolumn{1}{l|}{0}      & \multicolumn{1}{l|}{0.04}    & \multicolumn{1}{l|}{0}      & \multicolumn{1}{l|}{0}       & \multicolumn{1}{l|}{0}      & \multicolumn{1}{l|}{0}      & \multicolumn{1}{l|}{0}       & \multicolumn{1}{l|}{-0.02}   &  &  &  &  &  &  &  &  &  &  &  &  &  &  &  &  &  &  &  &  &  &  &  &  &  &  &  &  &  &  &  &  &  &  &  &  &  &  &  &  &  &  &  &  \\ \cline{1-11}
                               &                             &                              &                             &                              &                             &                              &                             &                             &                              &                              &  &  &  &  &  &  &  &  &  &  &  &  &  &  &  &  &  &  &  &  &  &  &  &  &  &  &  &  &  &  &  &  &  &  &  &  &  &  &  &  &  &  &  &  \\ \cline{1-11}
\multicolumn{1}{|l|}{Instance} & \multicolumn{2}{c|}{540}                                   & \multicolumn{2}{c|}{541}                                   & \multicolumn{2}{c|}{542}                                   & \multicolumn{2}{c|}{350}                                  & \multicolumn{2}{c|}{351}                                    &  &  &  &  &  &  &  &  &  &  &  &  &  &  &  &  &  &  &  &  &  &  &  &  &  &  &  &  &  &  &  &  &  &  &  &  &  &  &  &  &  &  &  &  \\ \cline{1-11}
\multicolumn{1}{|l|}{Method}     & \multicolumn{1}{l|}{CoBiT}  & \multicolumn{1}{l|}{GB}      & \multicolumn{1}{l|}{CoBiT}  & \multicolumn{1}{l|}{GB}      & \multicolumn{1}{l|}{CoBiT}  & \multicolumn{1}{l|}{GB}      & \multicolumn{1}{l|}{CoBiT}  & \multicolumn{1}{l|}{GB}     & \multicolumn{1}{l|}{CoBiT}   & \multicolumn{1}{l|}{GB}      &  &  &  &  &  &  &  &  &  &  &  &  &  &  &  &  &  &  &  &  &  &  &  &  &  &  &  &  &  &  &  &  &  &  &  &  &  &  &  &  &  &  &  &  \\ \cline{1-11}
\multicolumn{1}{|l|}{Time (S)}     & \multicolumn{1}{l|}{1020.46} & \multicolumn{1}{l|}{2862.82} & \multicolumn{1}{l|}{325.27} & \multicolumn{1}{l|}{3303.98} & \multicolumn{1}{l|}{371.05} & \multicolumn{1}{l|}{1101.4}  & \multicolumn{1}{l|}{8.58}   & \multicolumn{1}{l|}{7200}   & \multicolumn{1}{l|}{15.58}   & \multicolumn{1}{l|}{7200}    &  &  &  &  &  &  &  &  &  &  &  &  &  &  &  &  &  &  &  &  &  &  &  &  &  &  &  &  &  &  &  &  &  &  &  &  &  &  &  &  &  &  &  &  \\ \cline{1-11}
\multicolumn{1}{|l|}{Gap (\%)}    & \multicolumn{1}{l|}{0}      & \multicolumn{1}{l|}{0}       & \multicolumn{1}{l|}{0}      & \multicolumn{1}{l|}{0}       & \multicolumn{1}{l|}{0}      & \multicolumn{1}{l|}{0}       & \multicolumn{1}{l|}{0}      & \multicolumn{1}{l|}{1.81}   & \multicolumn{1}{l|}{0}       & \multicolumn{1}{l|}{14.18}   &  &  &  &  &  &  &  &  &  &  &  &  &  &  &  &  &  &  &  &  &  &  &  &  &  &  &  &  &  &  &  &  &  &  &  &  &  &  &  &  &  &  &  &  \\ \cline{1-11}
                               &                             &                              &                             &                              &                             &                              &                             &                             &                              &                              &  &  &  &  &  &  &  &  &  &  &  &  &  &  &  &  &  &  &  &  &  &  &  &  &  &  &  &  &  &  &  &  &  &  &  &  &  &  &  &  &  &  &  &  \\ \cline{1-11}
\multicolumn{1}{|l|}{Instance} & \multicolumn{2}{c|}{352}                                   & \multicolumn{2}{c|}{450}                                   & \multicolumn{2}{c|}{451}                                   & \multicolumn{2}{c|}{452}                                  & \multicolumn{2}{c|}{550}                                    &  &  &  &  &  &  &  &  &  &  &  &  &  &  &  &  &  &  &  &  &  &  &  &  &  &  &  &  &  &  &  &  &  &  &  &  &  &  &  &  &  &  &  &  \\ \cline{1-11}
\multicolumn{1}{|l|}{Method}     & \multicolumn{1}{l|}{CoBiT}  & \multicolumn{1}{l|}{GB}      & \multicolumn{1}{l|}{CoBiT}  & \multicolumn{1}{l|}{GB}      & \multicolumn{1}{l|}{CoBiT}  & \multicolumn{1}{l|}{GB}      & \multicolumn{1}{l|}{CoBiT}  & \multicolumn{1}{l|}{GB}     & \multicolumn{1}{l|}{CoBiT}   & \multicolumn{1}{l|}{GB}      &  &  &  &  &  &  &  &  &  &  &  &  &  &  &  &  &  &  &  &  &  &  &  &  &  &  &  &  &  &  &  &  &  &  &  &  &  &  &  &  &  &  &  &  \\ \cline{1-11}
\multicolumn{1}{|l|}{Time (S)}     & \multicolumn{1}{l|}{66.00}  & \multicolumn{1}{l|}{7200}    & \multicolumn{1}{l|}{80.87}  & \multicolumn{1}{l|}{7200}    & \multicolumn{1}{l|}{25.06}  & \multicolumn{1}{l|}{7200}    & \multicolumn{1}{l|}{93.06} & \multicolumn{1}{l|}{7200}   & \multicolumn{1}{l|}{1612.18} & \multicolumn{1}{l|}{7200}    &  &  &  &  &  &  &  &  &  &  &  &  &  &  &  &  &  &  &  &  &  &  &  &  &  &  &  &  &  &  &  &  &  &  &  &  &  &  &  &  &  &  &  &  \\ \cline{1-11}
\multicolumn{1}{|l|}{Gap (\%)}    & \multicolumn{1}{l|}{0}      & \multicolumn{1}{l|}{2.4}     & \multicolumn{1}{l|}{0}      & \multicolumn{1}{l|}{3.13}    & \multicolumn{1}{l|}{0}      & \multicolumn{1}{l|}{0}       & \multicolumn{1}{l|}{0}      & \multicolumn{1}{l|}{6.56}   & \multicolumn{1}{l|}{0}       & \multicolumn{1}{l|}{3.7}     &  &  &  &  &  &  &  &  &  &  &  &  &  &  &  &  &  &  &  &  &  &  &  &  &  &  &  &  &  &  &  &  &  &  &  &  &  &  &  &  &  &  &  &  \\ \cline{1-11}
                               &                             &                              &                             &                              &                             &                              &                             &                             &                              &                              &  &  &  &  &  &  &  &  &  &  &  &  &  &  &  &  &  &  &  &  &  &  &  &  &  &  &  &  &  &  &  &  &  &  &  &  &  &  &  &  &  &  &  &  \\ \cline{1-5}
\multicolumn{1}{|l|}{Instance} & \multicolumn{2}{c|}{551}                                   & \multicolumn{2}{c|}{552}                                   &                             &                              &                             &                             &                              &                              &  &  &  &  &  &  &  &  &  &  &  &  &  &  &  &  &  &  &  &  &  &  &  &  &  &  &  &  &  &  &  &  &  &  &  &  &  &  &  &  &  &  &  &  \\ \cline{1-5}
\multicolumn{1}{|l|}{Method}     & \multicolumn{1}{l|}{CoBiT}  & \multicolumn{1}{l|}{GB}      & \multicolumn{1}{l|}{CoBiT}  & \multicolumn{1}{l|}{GB}      &                             &                              &                             &                             &                              &                              &  &  &  &  &  &  &  &  &  &  &  &  &  &  &  &  &  &  &  &  &  &  &  &  &  &  &  &  &  &  &  &  &  &  &  &  &  &  &  &  &  &  &  &  \\ \cline{1-5}
\multicolumn{1}{|l|}{Time (S)}     & \multicolumn{1}{l|}{77.88}  & \multicolumn{1}{l|}{7200}    & \multicolumn{1}{l|}{15.79}  & \multicolumn{1}{l|}{7200}    &                             &                              &                             &                             &                              &                              &  &  &  &  &  &  &  &  &  &  &  &  &  &  &  &  &  &  &  &  &  &  &  &  &  &  &  &  &  &  &  &  &  &  &  &  &  &  &  &  &  &  &  &  \\ \cline{1-5}
\multicolumn{1}{|l|}{Gap (\%)}    & \multicolumn{1}{l|}{0}      & \multicolumn{1}{l|}{2.46}    & \multicolumn{1}{l|}{0}      & \multicolumn{1}{l|}{7.81}    &                             &                              &                             &                             &                              &                              &  &  &  &  &  &  &  &  &  &  &  &  &  &  &  &  &  &  &  &  &  &  &  &  &  &  &  &  &  &  &  &  &  &  &  &  &  &  &  &  &  &  &  &  \\ \cline{1-5}
\end{tabular}}
\caption{\ahmadreza{Comparison between the method proposed by \cite{GeerBranch}, denoted by GB, and CoBiT.}}\label{Table:comparison of CoBiT with GB}
\end{table}

In this section, we compare \ahmadreza{CoBiT} with the algorithm developed by \cite{GeerBranch}. We follow the same steps in generating random instances. More specifically, we consider $I,N\in\{3,4,5\}$. We consider $C = \mathcal{I}\cup \{0\}$, where $0$ refers to the non-purchase alternative. Given the number of alternatives $I$ and the number of customers $N$, we randomly generate three instances from different seeds. For each instance, $q_{in} \sim U(-7.0,7.0)$ and $\beta^p_{i}\sim U(-0.01,-0.001)$, where $U(a,b)$ denotes a uniform distribution in the range $(a,b).$ We have put a time limit of $7200$ seconds on both methods. Regarding the quality, we set the relative error to be $\epsilon=0.00001$ in both algorithms.}

\ahmadreza{The results obtained by applying the algorithm proposed by Geer and Branch \cite{GeerBranch}, referred to as GB, as well as applying \ahmadreza{CoBiT}, are summarized in Table \ref{Table:comparison of CoBiT with GB}. The table shows a comparison of the performance of the two algorithms in different instances. The `Instance' rows provide the name of the instance, which is composed of three numbers. The first two numbers represent $I$ and $N$, respectively, and the last number indicates the index of the random instance in this class of instances, which can take values in the set $\{0,1,2\}$.  To compare the solutions, we report the time (in seconds) taken by the algorithm to solve the instance. In case the time limit of $7200$ seconds is reached, we report the optimality gap, calculated by 
$
\frac{U-L}{L},
$
where $U$ and $L$ denote the upper and lower bounds, respectively. 

From this table, interesting observations can be made: (i) GB reports negative optimality gaps in 4 instances, meaning the upper bounds are lower than the lower bounds for these instances. Such errors occur in branching techniques due to numerical rounding. Therefore, we can safely assume that GB solves these instances to optimality. (ii) \ahmadreza{CoBiT} can solve all the instances to optimality within the given time limit, but GB can only solve 16 out of 27 instances. Also, among the instances solved by both algorithms, \ahmadreza{CoBiT} has a significantly lower time than GB, except the instance 531. (iii) The performance of GB varies significantly across different classes of instances. While GB may appear to perform well on instances with $N=3$ based on the experiments conducted by \cite{GeerBranch}, it is important to note that the complexity of GB is $\mathcal{O}(\epsilon ^{-N} I^{5.5+3N})$. The authors of the study used $\epsilon=0.01$ in their numerical experiments, while we consider $\epsilon=0.00001$ to ensure that the obtained solutions are indeed optimal. (iv) As expected, \ahmadreza{CoBiT} is quite sensitive to an increase in $I$ rather than an increase in $N$. We see that for $I=3,\; 4,$ and 5, the instances are solved in less than $67$, $95$, and $1615$ seconds, respectively, while changes in $N$ do not necessarily change the solution time.  
 }

\subsection{\ahmadreza{CoBiT} versus restrictive pricing problems} \label{Sec:parking choice model}
In this section, we investigate \ahmadreza{CoBiT} from a problem definition perspective by comparing pricing problems that differ in terms of the attributes and price sensitivity parameters included in the utility functions. 

\subsubsection{The parking choice model}

The selection of this case study is motivated by the availability of a published, non-trivial, disaggregate parking choice model by \cite{ibeas2014modelling}, \ahmadreza{which} we can use to characterize the demand. Furthermore, this case study has been recently used by \cite{paneque2018lagrangian} to demonstrate how to integrate advanced discrete choice models in pricing problems using a mixed integer linear programming (MILP) formulation. 
The parking choice consists \ahmadreza{of} three services: (1) paid on-street parking (PSP), (2) paid parking in an underground car park (PUP), and (3) free on-street parking (FSP). The latter does not provide any revenue to the operator. Table \ref{explanatoryvars} shows all explanatory variables used in the utility functions of the logit model. These are features related to the age of the vehicle, the income of customers, the type of trip, the access time to the destination from the parking, and information \ahmadreza{on} whether the customer is a resident or not. 

\small
\begin{table}[ht]
    \centering
    \begin{tabular}{l|p{13.5cm}}
       Features     & Definition \\
        \hline
        $ASC_{PSP}$ & Alternative specific constant for the PSP alternative. \\ 
        \hdashline
        $ASC_{PUP}$ & Alternative specific constant for the PUP alternative. \\ 
        \hdashline
        $AT_{FSP}$ & Access time to the free on-street parking. \\ 
        \hdashline
        $AT_{PSP}$ & Access time to the paid on-street parking. \\ 
        \hdashline
        $AT_{PUP}$ & Access time to the paid underground parking. \\
        \hdashline
        $TD_{FSP}$ & Access time to the destination from the free on-street parking. \\ 
        \hdashline
        $TD_{PSP}$ & Access time to the destination from the paid on-street parking. \\ 
        \hdashline
        $TD_{PUP}$ & Access time to the destination from the paid underground parking. \\ 
        \hdashline
        $Origin$ & Dummy parameter that is 1 if the origin of the trip is internal
to the town. \\
        \hdashline
        $\textbf{p}_{\textbf{PSP}}$ & \textbf{Fee for the paid on-street parking.}  \\ 
        \hdashline
        $\textbf{p}_{\textbf{PUP}}$ & \textbf{Fee for the paid underground parking.} \\ 
        \hdashline
        $LowInc$ & Dummy parameter that is 1 if the income of the customer is below
1200\euro{}/month. \\ 
\hdashline
$Residence$ & Dummy parameter that is 1 if the customer is a resident.\\ 
\hdashline
$AgeVeh_{\leq 3}$ & Dummy parameter that is 1 if the age of the vehicle is lower than 3 years. \\
\hline
\end{tabular}
    \caption{Features used in the parking choice model.}
    \label{explanatoryvars}
\end{table}
\normalsize

Following the logit model proposed by \cite{ibeas2014modelling}, we build the following three utility specifications: 
\begin{align*}
    V_{FSP,n} &= \beta^p_{FSP,n} \times p_{FSP} + q_{FSP,n}  = q_{FSP,n} \\
    V_{PSP,n} &= \beta^p_{PSP,n} \times p_{PSP} + q_{PSP,n}, \\
    V_{PUP,n} &= \beta^p_{PUP,n} \times p_{PUP} + q_{PUP,n}.
\end{align*}
The utility specification of the free on-street parking only contains the exogenous part $q_{FSP,n}$ since there is no fee to pay for that option ($p_{FSP}=0$). 
The price sensitivities parameters $\beta^p_{PSP,n}$ and $\beta^p_{PUP,n}$ are then further expressed as:
\begin{align}
\label{betaprices}
\beta^p_{PSP,n} &=\beta_{FEE}+\beta_{FEE_{PSP(LowInc)}}\times LowInc_n + \beta_{FEE_{PSP(Resident)}}\times Residence_n\\
\label{betaprices2}
\beta^p_{PUP,n}&=\beta_{FEE}+\beta_{FEE_{PUP(LowInc)}} \times LowInc_n+ \beta_{FEE_{PUP(Resident)}} \times Residence_n.
\end{align}
The exogenous parts of utilities are modeled as:
\begin{align*}
q_{FSP,n}&= &\beta_{AT} \times AT_{FSP} + &\beta_{TD} \times TD_{FSP} &&+\beta_{Origin}\times Origin_n, \\
q_{PSP,n}&= ASC_{PSP} +& \beta_{AT}\times AT_{PSP} + &\beta_{TD}\times TD_{TSP},&&\\
q_{PUP,n}&= ASC_{PUP} +& \beta_{AT}\times AT_{PUP} +& \beta_{TD}\times TD_{PUP} && + \beta_{AgeVeh_{\leq 3}}\times AgeVeh_{{\leq 3}_n}.
\end{align*}
The values of coefficient parameters used in \cite{ibeas2014modelling} are depicted in Table \ref{params}. Note that in \cite{ibeas2014modelling}, $\beta_{AT}$ and $\beta_{FEE}$ are assumed to be normally distributed and correlated, with $cov(\beta_{AT}; \beta_{FEE}) = -12.8.$

\small
\begin{table}[ht]
    \centering
    \begin{tabular}{c|c}
          & Mixed Logit  \\
         \hline
    $ASC_{PSP}$  & 32  \\
    $ASC_{PUP}$  & 34  \\
    $\beta_{TD}$    & -0.612 \\
    $\beta_{Origin}$  & -5.762  \\
    $\left[\begin{matrix}\beta_{AT} \\ \beta_{FEE}\end{matrix}\right]$  & $ \sim Normal \left(\left[\begin{matrix}-0.788 \\ -32.3 \end{matrix} \right], \left[\begin{matrix} 1.1236 & -12.8\\-12.8 & 201.64 \end{matrix} \right]\right)$  \\
    $\beta_{FEE_{PSP(LowInc)}}$  & -10.995  \\
        $\beta_{FEE_{PSP(Resident)}}$  & -11.44  \\
            $\beta_{FEE_{PUP(LowInc)}}$  & -13.729   \\
        $\beta_{FEE_{PUP(Resident)}}$  & -10.668   \\
   $\beta_{AgeVeh_{\leq 3}}$  & 4.037  \\
\end{tabular}
    \caption{Values of coefficient parameters. }\label{Tab:parameters of PArking}
    \label{params}
\end{table}

\normalsize

The pricing problem is to determine the optimal prices (or parking fees) of the two paid parking services, \textit{i.e.}, $p_{_{PSP}}$ and $p_{_{PUP}}$, so that the revenue of the operator is maximized. Since the purpose is to show the practicality of \ahmadreza{CoBiT}, we consider an unlimited capacity for the parking services. In the pricing problem, $p_{_{PSP}}$ and $p_{_{PUP}}$ are the only endogenous variables, and all others are exogenous demand variables for which values are given.

\subsubsection{Numerical results}\label{Sec:numerical discrete}
{As mentioned, \ahmadreza{CoBiT} is capable of solving pricing problems under discrete mixed logit models. Unlike existing contributions, \ahmadreza{CoBiT} is capable to handle heterogeneous price sensitivity parameters in the utility function, better reflecting the demand models that have been used in the DCM literature.} In this section, we perform experiments with two main goals in mind:
\begin{itemize}
\item[\textbf{G1}] First, we want to show that \ahmadreza{CoBiT} outperforms the global optimizers SCIP 7.0, BARON 21.1.13, Couenne, and ANTIGONE 37.1, both in time and optimality gap. To do so, we assume an MNL model, with $\beta_{AT}$ and $\beta_{FEE}$ fixed to their mean values (\textit{i.e.}, $\beta_{AT}=-0.788$ and $\beta_{FEE}=-32.3$) and we use instances with 10 and 50 customers. It is worth noting that under these assumptions, the price sensitivity parameters ($\beta^p_{PSP,n}$ and $\beta^p_{PUP,n}$) are still both product and customers dependent (see Equations (\ref{betaprices})-(\ref{betaprices2})), which means that state-of-the-art methods cannot be used to solve these problems.

\item[\textbf{G2}] Second, we show the consequences, mainly in terms of lost revenues, that would arise from optimizing under simplified assumptions regarding the price sensitivity parameters, either by assuming a single value instead of a distribution or by neglecting that the price sensitivity parameters can be both product and customers dependent. 
\end{itemize}

\begin{table}[h!]
\centering
\begin{tabular}{|c|c|c|c|c|c|c|c|}
\hline
$\text{Number of}\atop\text{customers}$ &Method& $\mbox{Best}\atop \mbox{solution}$ & Upper bound & Opt. Gap&$p^*_{PSP}$&$p^*_{PUP} $&$\mbox{Time}\atop \mbox{(Seconds)}$\\\hline
 10&  \ahmadreza{CoBiT}   &  6.36 & 6.36& 0.00\% &0.7036&0.7137& \ahmadreza{514.45}\\
10& SCIP     & 0.0 & $3\times 10^{19}$ & $10^{21} $\%&0&0&7,200\\
10& BARON    & 6.01 &6.47& 7.65\%&0.6117&0.6202&7,200\\
10& Couenne &6.36&774.91& $1.2\times 10^4\%$&0.7025&0.7176&7,200\\
10& ANTIGONE&6.35$^*$&-$^*$& $0.00\%$ &0.6551&0.6638&588\\\hline
 50&  \ahmadreza{CoBiT}   &  31.93 & 31.93& 0.00\% &0.7142&0.7210& \ahmadreza{780.02}\\
50& SCIP     & 31.74 & $ 10^{20}$ & $10^{21} $\%&0.6531&0.6642&7,200\\
50& BARON    & 30.51 &587,986& $1.9\times 10^6$\%&0.6178&0.6291&7,200\\
50& Couenne&31.93&476,747.02& $1.5\times 10^6\%$&0.7142&0.7210&7,200\\
50& ANTIGONE&31.74&10,000& $3.1\times 10^4\%$ &0.6531&0.6642&7,200\\\hline
\end{tabular}
\caption{Information obtained on solving \eqref{Eq:MixedMNL} with degenerate mixing probability measures  with $N=10$ and $N=50$ customers.\\
$^*$:Looking at the log of ANTIGONE, it seems that there is a bug in the solver. After 525 Seconds, the solver reports a solution whose objective value is 7.449, and after 588 Seconds ANTIGONE proves optimality of that solution. We have reported this bug.}\label{Tabel:comparison for MNL N=10}
\end{table}

\textbf{G1.} Table \ref{Tabel:comparison for MNL N=10} provides the comparison between \ahmadreza{CoBiT} and the other solvers. Since the problem contains a fraction and exponential functions, its solution is quite sensitive to the errors. That made the solution obtained by ANTIGONE unreliable. As one can see, \ahmadreza{CoBiT} is the only solver that can solve the problem to optimality. Couenne solver can find the optimal solution within two hours but it has difficulty proving optimality of it. On the other hand, BARON is the solver with second best optimality gap. One of the main reasons why solvers are unable to solve these instances is the way they construct an upper bound. Most of the solvers use a technique called \textit{term-based underestimates}, where they introduce new variables to represent each nonlinear terms (in our instances, the nonlinear terms are the fraction and the exponential functions). Using this technique and then convexification of the problem results in loose relaxations (convexification of each term comes with some gap; hence putting all the relaxations into one problem aggregate the error). However, in \ahmadreza{CoBiT}, we only consider the summation of the exponential functions as one nonlinear term and the fraction as the other one. Therefore, our convexification is tighter. 

To have a better understanding on how \ahmadreza{CoBiT} converges to the optimal solution, we illustrate how the optimality gap is reduced over time. In all instances, \ahmadreza{CoBiT} finds optimal solutions in the first few iterations and attempts to close the optimality gap in the later ones. 
As expected, \ahmadreza{CoBiT} converges faster for the instance with $N=10$ customer classes, because in this instance each iteration can be solved much faster. 
\begin{figure}[t]
\centering
%
%
\begin{tikzpicture}
\footnotesize
\begin{axis}[%
width=12cm,
height=3cm,
at={(1.55cm,0.434cm)},
scale only axis,
unbounded coords=jump,
xmin=0,
xmax=800,
xlabel style={font=\color{white!15!black}},
xlabel={Time (S)},
ymin=0,
ymax=100,
ylabel style={font=\color{white!15!black}},
ylabel={Optimality Gap (\%)},
axis background/.style={fill=white},
legend style={legend cell align=left, align=left, draw=white!15!black}
]

\addplot [color=red, solid, line width=2.0pt]
  table[row sep=crcr]{%
20.78	94.5151155879402\\
22.483	88.9915900337971\\
24.499	83.9330346616364\\
27.389	68.2779218737719\\
31.749	60.3678377741099\\
37.326	55.089680106893\\
43.92	46.0737247504519\\
54.545	32.1000831563575\\
71.375	23.2617727552106\\
97.983	14.0743314087379\\
125.357	9.3931510402339\\
165.248	4.59086228985074\\
208.186	2.50174878763822\\
267.622	0.996769655194085\\
319.997	0.384025905636201\\
380.388	0.384025905636201\\
427.464	0.384025905636201\\
460.746	0.221172512988177\\
480.215	0.160023893547957\\
494.465	0.106420605041216\\
506.247	0.103748300335606\\
514.451	0\\
};
\addlegendentry{$\text{LiBiT}\text{ (N=10)}$}


\addplot [color=blue, solid, line width=2.0pt, mark=star, mark options={solid, blue}]
  table[row sep=crcr]{%
 1.985	94.3388010325556\\
6.376	87.4896619132353\\
10.408	82.185949481281\\
15.814	64.9739387560768\\
22.876	54.3987746203578\\
33.282	48.2834661454418\\
46.955	39.4104270034581\\
69.549	27.560275868981\\
95.642	19.4472632515253\\
136.313	11.4505849006374\\
180.767	7.40997099778882\\
253.282	3.50250245861078\\
332.846	1.91270522352999\\
450.521	0.566263475379445\\
560.443	0.566263475379445\\
661.617	0.564071083605293\\
733.008	0.310693234278984\\
778.039	0.0460402272570713\\
780.023	0\\
};
\addlegendentry{$\text{LiBiT}\text{ (N=50)}$}

\end{axis}
\end{tikzpicture}%
    \caption{\ahmadreza{The convergence behaviour of \ahmadreza{CoBiT} over time on the parking case with $N=10$ and $N=50$ customers.}}
    \label{fig:convergence n=10}
\end{figure}


\begin{figure}[t]
\centering
%
%
\begin{tikzpicture}
\footnotesize
\begin{axis}[%
width=12cm,
height=3cm,
at={(1.55cm,0.434cm)},
scale only axis,
xmin=0,
xmax=24,
xlabel style={font=\color{white!15!black}},
xlabel={Iteration Number},
ymin=0,
ymax=250,
ylabel style={font=\color{white!15!black}},
ylabel={Number of Nodes},
axis background/.style={fill=white},
legend style={at={(3.25cm,2.7cm)},legend cell align=left, align=left, draw=white!15!black}
]
\addplot [color=red, solid, line width=2.0pt]
  table[row sep=crcr]{%
  1    2  \\
  2    4  \\
  3    4  \\
  4    6  \\
  5   10  \\
  6   16  \\
  7   16  \\
  8   28  \\
  9   40  \\
 10   68  \\
 11   72  \\
 12  116  \\
 13  124  \\
 14  188  \\
 15  186  \\
 16  224  \\
 17  186  \\
 18  168  \\
 19   98  \\
 20   70  \\
 21   58  \\
 22   42  \\
 23 1\\
};
\addlegendentry{$\text{LiBiT} \text{ (N=10)}$}


\addplot [color=blue, solid, line width=2.0pt, mark=star, mark options={solid, blue}]
  table[row sep=crcr]{%
  1    2\\
  2    4\\
  3    4\\
  4    6\\
  5    8\\
  6   12\\
  7   16\\
  8   26\\
  9   30\\
 10   46\\
 11   50\\
 12   80\\
 13   88\\
 14  126\\
 15  118\\
 16  108\\
 17   76\\
 18   48\\
 19    2\\
};
\addlegendentry{$\text{LiBiT} \text{ (N=50)}$}


\end{axis}
\end{tikzpicture}%
    \caption{\ahmadreza{The explored  number of nodes in each iteration of CoBiT for the parking case with $N=10$ and $N=50$ customers. }}\label{Fig:explored nodes}
    \end{figure}
Next to the optimality gap, it is also interesting to see the number of nodes generated in each iteration of \ahmadreza{CoBiT}. Figure \ref{Fig:explored nodes} provides this information. As one can see, the behavior of \ahmadreza{CoBiT} in both instances is similar. More specifically, we see that the number of nodes that are explored has its peak around the 15th iteration, and then it decreases. Moreover, an interesting observation is on the solid blue curve (corresponding to 50 customers), which is below the solid red curve (corresponding to 10 customers). The reason is the shape of the expected revenue function. As we can see in Appendix \ref{Appendix: figure for objective functions Mixed}, the expected revenue function is rather flat when considering the instance with 10 customers compared to the one with 50 customers. Therefore, \ahmadreza{CoBiT} needs small partitions to make sure where the optimal solution is in the instance with 10 customers.


\textbf{G2.} In this part, we focus on evaluating the impact of simplified assumptions on the quality of the solution. One possible way to simplify the pricing model is by assuming that the price sensitivity is customer-independent (as suggested in \cite{GeerBranch}). Mathematically speaking, this simplification reduces the dimensionality of $\beta^p$ from $\mathbb{R}^{C\times N}$ to $\mathbb{R}^C$. To investigate the effect of this simplification, we restrict ourselves to three columns of $\beta^p\in \mathbb{R}^{C\times N}$, which we refer to by Classes 1, 2, and 3.

To assess the impact of this simplification, we evaluate the optimal prices obtained by restricting ourselves to these three classes on the revenue function of the parking case study (where $\beta^p\in \mathbb{R}^{C\times N}$). Table \ref{Tab:GB assumption} displays the results obtained from this analysis. As the table illustrates, the optimal prices obtained are not optimal in all cases and can be far from optimal in some cases. Thus, neglecting customer-specific price sensitivity in modeling the problem can result in solutions that may not even be locally optimal, indicating the importance of considering customer-specific price sensitivity in the problem's modeling.
\begin{table}[h!]
\centering
\ahmadreza{
\begin{tabular}{|c|c|c|c|c|c|}
\hline
Customers' price	& 	& 	 & Evaluation for problem  &	 \\
 sensitivity	& $P_PSP$	& $P_PUP$	 & with 10 classes &	Gap(\%) \\
\hline
Class 1	&0.62 &	0.62	&6.20 &	2.60 \\
\hline
Class 2	&0.77 &	0.77	&5.38	 &15.43 \\
\hline
Class 3	&0.82	& 0.82	&3.73	 &41.40 \\
\hline
\end{tabular}}
\caption{\ahmadreza{The result of ignoring the customers' sensitivity by only considering one customer class to calculate $\beta^p$.} }\label{Tab:GB assumption}
\end{table}
 

As mentioned, the original distributions of $\beta_{AT}$ and $\beta_{FEE}$ are continuous, while the above-mentioned results are achieved by considering a single value (mean value) of these parameters. Thanks to \ahmadreza{CoBiT}, we are able to integrate a discrete mixed logit model within our pricing problem, i.e., to discretize the distribution and approximate the integral. 

Let us consider the situation with $N=10$ customer classes. We limit the support set of $\left[\beta_{AT} \atop \beta_{FEE}\right]$ to its 0.99 confidence set, i.e., $[-3.6,1.94]\times [-68.52,3.92].$ To discretize the distribution and approximate the integral, we break the length and the width of the confidence set into $n$ parts with the same length; hence $n^2$ break points. We check the solution obtained by this approximation when $n^2$ varies in $\{9,16,25,49,64, 100, 121, 144, 169, \ahmadreza{400, 900}\}$. As mentioned before, the performance of the available solvers \ahmadreza{depends} on the number of nonlinear terms in the optimization problem. Discretization of the continuous distribution increases the number of nonlinear terms dramatically, hence negatively \ahmadreza{affecting} the performance of the solvers. Since we have seen that the available solvers cannot tackle the simple MNL model, we only apply \ahmadreza{CoBiT} for the discrete cases.

Table \ref{Tab:sensitivity to samples} shows how the optimal value and optimal solutions change when the number of breakpoints increases. In the fifth column, we also report the value of the objective function of the continuous mixed logit model for the obtained optimal prices.

\begin{table}[h]
\centering
\begin{tabular}{|c|c|c|c|c|c|}
\hline
Num.   break points & Opt. value  &$p^*_{PSP}$   & $p^*_{PUP}$  & $\mbox{Expected revenue}\atop \mbox{of continuous ML}$ &Time (Seconds) \\ \hline
1 (MNL case) &  6.36 & 0.70&0.71& 4.43& \ahmadreza{514.45}\\\hline
  9             & \ahmadreza{6.80}     & \ahmadreza{0.53}    & \ahmadreza{0.75}     &  \ahmadreza{5.05}   & \ahmadreza{62.16}\\ \hline
16             & \ahmadreza{5.32} &  \ahmadreza{0.49} & \ahmadreza{0.64}& \ahmadreza{5.07}   & \ahmadreza{184.42} \\ \hline
25            & \ahmadreza{5.11}  &  \ahmadreza{0.56} & \ahmadreza{0.70}& \ahmadreza{5.00}   &\ahmadreza{225.81} \\ \hline 
49            & \ahmadreza{5.09}    & \ahmadreza{0.51} & \ahmadreza{0.66} &  \ahmadreza{5.09} & \ahmadreza{323.58} \\ \hline 
64            & \ahmadreza{5.09}    & \ahmadreza{0.50} & \ahmadreza{0.66}&  \ahmadreza{5.08} & \ahmadreza{362.45} \\ \hline
100            & \ahmadreza{5.06}    & \ahmadreza{0.50} & \ahmadreza{0.66}&  \ahmadreza{5.08} & \ahmadreza{278.66}\\ \hline
121            &\ahmadreza{5.06}    & \ahmadreza{0.50} & \ahmadreza{0.67}&  \ahmadreza{5.08} & \ahmadreza{53.52}\\ \hline
144            & \ahmadreza{5.07}    & \ahmadreza{0.50} & \ahmadreza{0.66}&  \ahmadreza{5.08} & \ahmadreza{407.09}\\ \hline
169           & \ahmadreza{5.07}    & \ahmadreza{0.50} & \ahmadreza{0.66}&  \ahmadreza{5.08} & \ahmadreza{84.13}\\ \hline
\ahmadreza{400}           & \ahmadreza{5.08}    & \ahmadreza{0.50} & \ahmadreza{0.67}&  \ahmadreza{5.08} & \ahmadreza{989.66}\\ \hline
\ahmadreza{900}          & \ahmadreza{5.08}    & \ahmadreza{0.50} & \ahmadreza{0.67}&  \ahmadreza{5.08} & \ahmadreza{2,912.58}\\ \hline
\end{tabular}
\caption{Sensitivity of the optimal prices to the number of break points.}\label{Tab:sensitivity to samples}
\end{table}

We see that the sequence of solutions converges as we increase the number of breakpoints. Doing so, we better approximate the continuous mixed logit model, which is the best choice model reported in \cite{ibeas2014modelling}. Accordingly, our pricing model better reflects the heterogeneity of the parking users' behaviors. While we expect this additional complexity in the demand model to come with an increased computational time, we can see that this is not always true. \ahmadreza{In the end}, it depends, as mentioned already, on the resulting shape of the expected revenue function.

In Table \ref{tab:solutions of different models}, we show the revenue, as well as the markets shares obtained while using the optimal prices of the simple MNL, as well as several discrete mixed logit models, into the revenue maximization objective function of the continuous mixed logit demand. 

\begin{table}[h]
     \centering
     \begin{tabular}{|c|c|c|c|c|c|c|c|}
     \hline

          Model &$N$ & $p^*_{PSP}$& $p^*_{PUP}$ & Expected revenue & \multicolumn{3}{c|}{Market shares (\%)} \\
        &&&& &FSP&PSP&PUP\\\hline     \hline
           Simple MNL & 10& 0.70&0.71&4.43&36.35&1.41&60.66 \\\hline
           Discrete ML(9) & 10& \ahmadreza{0.53}    & \ahmadreza{0.75}     &  \ahmadreza{5.05}&\ahmadreza{14.20}&\ahmadreza{57.70}&\ahmadreza{26.53} \\\hline
           Discrete ML(16) & 10& \ahmadreza{0.49} & \ahmadreza{0.64}& \ahmadreza{5.07} & \ahmadreza{8.00}&\ahmadreza{48.75}&\ahmadreza{41.68} \\\hline
           Discrete ML(25) & 10& \ahmadreza{0.51} & \ahmadreza{0.66} &  \ahmadreza{5.09}&\ahmadreza{10.52}&\ahmadreza{47.50}&\ahmadreza{40.41} \\\hline
        Discrete ML(\ahmadreza{900}) & 10&  \ahmadreza{0.50} & \ahmadreza{0.67}&  \ahmadreza{5.08}&\ahmadreza{8.97}&\ahmadreza{52.63}&\ahmadreza{36.83} \\\hline
     \end{tabular}
     \caption{Expected revenue and market shares associated with the optimal prices of different logit models.}
    \label{tab:solutions of different models}
 \end{table}

On the parking choice instances of 10 customers, we see that assuming a simple MNL model leads to a significant drop in revenue (4.43 instead of 5.09). Even more interesting is to see that these assumptions regarding the demand also have an important impact on the market shares associated with the different parking options. We see that assuming a simple MNL model in our case would cause a significant shift in the market shares distribution. These are \ahmadreza{an} important consideration to take into account when deciding about pricing. Naturally, the magnitude of the revenue loss and market share shifts will depend on the realism of the assumptions regarding the demand model. With this work, we do not claim that a discrete mixed logit model should be used in all cases, but we offer an algorithm able to include this complex choice model if the analyst considers that this is the right choice model to use. Ultimately, the decision should be based on the trade-off between the demand model realism on one side,  and the complexity of the resulting pricing model on the other side, and this should be assessed on a case-by-case basis by the analyst. 

\color{black}
\section{Acknowledgment}
We would like to thank Dick den Hertog and Danique de Moor for discussing the convexification technique used in this paper.
\section{Conclusions}
Pricing problems under disaggregate demand assumptions is a challenging but relevant area of research due to its numerous applications. In this paper, we explored a static multi-product pricing problem under a discrete mixed logit model \ahmadreza{that can accommodate heterogeneous price sensitivity}. We designed an efficient iterative optimization algorithm that asymptotically converges to the optimal solution. We used linear optimization problems designed based on a trust-region approach to approximate the problem from below and therefore find a ``good" feasible solution. We then used \ahmadreza{convex} approximations as well as McCormick relaxation to obtain an upper bound on the optimal value of the nonlinear optimization problem.  Thanks to the branching method, we then tightened the optimality gap and proved \ahmadreza{the} asymptotic convergence of our algorithm.  

The effectiveness of this general algorithm was demonstrated on \ahmadreza{several} case studies. Benchmarks against solvers and existing contributions in the literature were performed. Our results showed that our algorithm \ahmadreza{could} find optimal prices, even for \ahmadreza{a} higher degree of segment asymmetry. Furthermore, we see that the considered  static multi-product pricing problems under a discrete mixed logit model can be considered as the hard instances for global optimization solvers and can be used as the test instances.

\ahmadreza{The computational complexity of the proposed algorithm is an important issue that needs further investigation. A formal computational complexity will show how the computation of our algorithm is linked to the number of customers as well as the number of alternatives.
Furthermore, this paper focused on the static pricing problem under discrete mixed logit demand. However, the proposed algorithm could be extended to other types of demand models and pricing problems. Exploring the generalization of the algorithm to other settings is an interesting direction for future research.
Finally, the partitioning method used in the proposed algorithm is a key component for ensuring the convergence of the algorithm. Investigating alternative partitioning methods and comparing their performance could be an interesting avenue for future research. }

\nolinenumbers
\bibliographystyle{informs2014}
\bibliography{references}

\begin{thebibliography}{53}
\providecommand{\natexlab}[1]{#1}
\providecommand{\url}[1]{\texttt{#1}}
\providecommand{\urlprefix}{URL }

\bibitem[{Ak{\c{c}}ay et~al.(2010)Ak{\c{c}}ay, Natarajan, \protect\BIBand{}
  Xu}]{akccay2010joint}
Ak{\c{c}}ay Y, Natarajan HP, Xu SH (2010) Joint dynamic pricing of multiple
  perishable products under consumer choice. \emph{Management Science}
  56(8):1345--1361.

\bibitem[{Akrotirianakis \protect\BIBand{}
  Floudas(2004)}]{akrotirianakis2004computational}
Akrotirianakis IG, Floudas CA (2004) Computational experience with a new class
  of convex underestimators: Box-constrained {NLP} problems. \emph{Journal of
  Global Optimization} 29(3):249--264.

\bibitem[{Aksoy-Pierson et~al.(2013)Aksoy-Pierson, Allon, \protect\BIBand{}
  Federgruen}]{aksoy2013price}
Aksoy-Pierson M, Allon G, Federgruen A (2013) Price competition under mixed
  multinomial logit demand functions. \emph{Management Science}
  59(8):1817--1835.

\bibitem[{Aydin \protect\BIBand{} Porteus(2008)}]{aydin2008joint}
Aydin G, Porteus EL (2008) Joint inventory and pricing decisions for an
  assortment. \emph{Operations Research} 56(5):1247--1255.

\bibitem[{Aydin \protect\BIBand{} Ryan(2000)}]{aydin2000product}
Aydin G, Ryan JK (2000) Product line selection and pricing under the
  multinomial logit choice model. \emph{Proceedings of the 2000 MSOM
  conference} (Citeseer).

\bibitem[{Belotti(2009)}]{belotti2009couenne}
Belotti P (2009) Couenne: a user’s manual. Technical report, Technical
  report, Lehigh University.

\bibitem[{Ben-Akiva \protect\BIBand{} Bierlaire(2003)}]{ben2003discrete}
Ben-Akiva M, Bierlaire M (2003) Discrete choice models with applications to
  departure time and route choice. \emph{Handbook of transportation science},
  7--37 (Springer).

\bibitem[{Bertsimas et~al.(2020)Bertsimas, Sian~Ng, \protect\BIBand{}
  Yan}]{bertsimas2020joint}
Bertsimas D, Sian~Ng Y, Yan J (2020) Joint frequency-setting and pricing
  optimization on multimodal transit networks at scale. \emph{Transportation
  Science} (Articles in Advance).

\bibitem[{Bezanson et~al.(2017)Bezanson, Edelman, Karpinski, \protect\BIBand{}
  Shah}]{Julia}
Bezanson J, Edelman A, Karpinski S, Shah V (2017) Julia: A fresh approach to
  numerical computing. \emph{SIAM Review} 59(1):65--98.

\bibitem[{Bortolomiol et~al.(2021)Bortolomiol, Lurkin, \protect\BIBand{}
  Bierlaire}]{bortolomiol2021simulation}
Bortolomiol S, Lurkin V, Bierlaire M (2021) A simulation-based heuristic to
  find approximate equilibria with disaggregate demand models.
  \emph{Transportation Science} 55(1).

\bibitem[{Boyd \protect\BIBand{} Vandenberghe(2004)}]{boyd2004convex}
Boyd S, Vandenberghe L (2004) \emph{Convex optimization} (Cambridge University
  Press).

\bibitem[{Conn et~al.(2000)Conn, Gould, \protect\BIBand{}
  Toint}]{conn2000trust}
Conn AR, Gould NI, Toint PL (2000) \emph{Trust region methods}, volume~1
  (SIAM).

\bibitem[{Daziano et~al.(2017)Daziano, Sarrias, \protect\BIBand{}
  Leard}]{daziano2017consumers}
Daziano RA, Sarrias M, Leard B (2017) Are consumers willing to pay to let cars
  drive for them? {A}nalyzing response to autonomous vehicles.
  \emph{Transportation Research Part C: Emerging Technologies} 78:150--164.

\bibitem[{Dong et~al.(2009)Dong, Kouvelis, \protect\BIBand{}
  Tian}]{dong2009dynamic}
Dong L, Kouvelis P, Tian Z (2009) Dynamic pricing and inventory control of
  substitute products. \emph{Manufacturing \& Service Operations Management}
  11(2):317--339.

\bibitem[{Du et~al.(2016)Du, Cooper, \protect\BIBand{} Wang}]{du2016optimal}
Du C, Cooper WL, Wang Z (2016) Optimal pricing for a multinomial logit choice
  model with network effects. \emph{Operations Research} 64(2):441--455.

\bibitem[{Dunning et~al.(2017)Dunning, Huchette, \protect\BIBand{}
  Lubin}]{JuMP}
Dunning I, Huchette J, Lubin M (2017) {JuMP}: A modeling language for
  mathematical optimization. \emph{SIAM Review} 59(2):295--320.

\bibitem[{Floudas et~al.(2005)Floudas, Akrotirianakis, Caratzoulas, Meyer,
  \protect\BIBand{} Kallrath}]{FLOUDAS20051185}
Floudas C, Akrotirianakis I, Caratzoulas S, Meyer C, Kallrath J (2005) Global
  optimization in the 21st century: Advances and challenges. \emph{Computers \&
  Chemical Engineering} 29(6):1185 -- 1202.

\bibitem[{Gallego \protect\BIBand{} Wang(2014)}]{gallego2014multiproduct}
Gallego G, Wang R (2014) Multiproduct price optimization and competition under
  the nested logit model with product-differentiated price sensitivities.
  \emph{Operations Research} 62(2):450--461.

\bibitem[{Gamrath et~al.(2020)Gamrath, Anderson, Bestuzheva, Chen, Eifler,
  Gasse, Gemander, Gleixner, Gottwald, Halbig et~al.}]{gamrath2020scip}
Gamrath G, Anderson D, Bestuzheva K, Chen WK, Eifler L, Gasse M, Gemander P,
  Gleixner A, Gottwald L, Halbig K, et~al. (2020) The scip optimization suite
  7.0 .

\bibitem[{Gilbert et~al.(2014)Gilbert, Marcotte, \protect\BIBand{}
  Savard}]{gilbert2014mixed}
Gilbert F, Marcotte P, Savard G (2014) Mixed-logit network pricing.
  \emph{Computational Optimization and Applications} 57(1):105--127.

\bibitem[{Hanson \protect\BIBand{} Martin(1996)}]{hanson1996}
Hanson W, Martin K (1996) Optimizing multinomial logit profit functions.
  \emph{Management Science} 42(7):992--1003.

\bibitem[{Hensher et~al.(2005)Hensher, Rose, \protect\BIBand{}
  Greene}]{hensher2005applied}
Hensher DA, Rose JM, Greene WH (2005) \emph{Applied choice analysis: a primer}
  (Cambridge University Press).

\bibitem[{Hess et~al.(2018)Hess, Daly, \protect\BIBand{}
  Batley}]{hess2018revisiting}
Hess S, Daly A, Batley R (2018) Revisiting consistency with random utility
  maximisation: theory and implications for practical work. \emph{Theory and
  Decision} 84(2):181--204.

\bibitem[{Higham(1999)}]{higham1999trust}
Higham DJ (1999) Trust region algorithms and timestep selection. \emph{SIAM
  Journal on Numerical Analysis} 37(1):194--210.

\bibitem[{Hopp \protect\BIBand{} Xu(2005)}]{hopp2005product}
Hopp WJ, Xu X (2005) Product line selection and pricing with modularity in
  design. \emph{Manufacturing \& Service Operations Management} 7(3):172--187.

\bibitem[{Huh \protect\BIBand{} Li(2015)}]{huh2015pricing}
Huh WT, Li H (2015) Pricing under the nested attraction model with a multistage
  choice structure. \emph{Operations Research} 63(4):840--850.

\bibitem[{Ibeas et~al.(2014)Ibeas, Dell’Olio, Bordagaray, \protect\BIBand{}
  Ort{\'u}zar}]{ibeas2014modelling}
Ibeas A, Dell’Olio L, Bordagaray M, Ort{\'u}zar JdD (2014) Modelling parking
  choices considering user heterogeneity. \emph{Transportation Research Part A:
  Policy and Practice} 70:41--49.

\bibitem[{Jalali et~al.(2019)Jalali, Carmen, Van~Nieuwenhuyse,
  \protect\BIBand{} Boute}]{jalali2019quality}
Jalali H, Carmen R, Van~Nieuwenhuyse I, Boute R (2019) Quality and pricing
  decisions in production/inventory systems. \emph{European Journal of
  Operational Research} 272(1):195--206.

\bibitem[{Li \protect\BIBand{} Huh(2011)}]{li2011pricing}
Li H, Huh WT (2011) Pricing multiple products with the multinomial logit and
  nested logit models: Concavity and implications. \emph{Manufacturing \&
  Service Operations Management} 13(4):549--563.

\bibitem[{Li \protect\BIBand{} Webster(2017)}]{li2017optimal}
Li H, Webster S (2017) Optimal pricing of correlated product options under the
  paired combinatorial logit model. \emph{Operations Research}
  65(5):1215--1230.

\bibitem[{Li et~al.(2019)Li, Webster, Mason, \protect\BIBand{}
  Kempf}]{li2019product}
Li H, Webster S, Mason N, Kempf K (2019) Product-line pricing under discrete
  mixed multinomial logit demand: Winner—2017 {M\&SOM} practice-based
  research competition. \emph{Manufacturing \& Service Operations Management}
  21(1):14--28.

\bibitem[{Li \protect\BIBand{} Kamargianni(2019)}]{li2019integrated}
Li W, Kamargianni M (2019) An integrated choice and latent variable model to
  explore the influence of attitudinal and perceptual factors on shared
  mobility choices and their value of time estimation. \emph{Transportation
  Science} Vol. 54, No. 1:62–83.

\bibitem[{Maddah \protect\BIBand{} Bish(2007)}]{maddah2007joint}
Maddah B, Bish EK (2007) Joint pricing, assortment, and inventory decisions for
  a retailer's product line. \emph{Naval Research Logistics (NRL)}
  54(3):315--330.

\bibitem[{Manski(1977)}]{manski1977structure}
Manski CF (1977) The structure of random utility models. \emph{Theory and
  Decision} 8(3):229--254.

\bibitem[{McCormick(1976)}]{mccormick1976computability}
McCormick GP (1976) Computability of global solutions to factorable nonconvex
  programs: Part i—convex underestimating problems. \emph{Mathematical
  Programming} 10(1):147--175.

\bibitem[{McFadden(1977)}]{mcfadden78}
McFadden D (1977) Modelling the choice of residential location. Cowles
  Foundation Discussion Papers 477, Cowles Foundation for Research in
  Economics, Yale University.

\bibitem[{McFadden \protect\BIBand{} Train(2000)}]{mcfadden2000mixed}
McFadden D, Train K (2000) Mixed {MNL} models for discrete response.
  \emph{Journal of Applied Econometrics} 15(5):447--470.

\bibitem[{McFadden \protect\BIBand{} Zarembka(1974)}]{mcfadden1974frontiers}
McFadden D, Zarembka P (1974) Frontiers in econometrics. \emph{Conditional
  Logit Analysis of Qualitative Choice Behavior} 105--142.

\bibitem[{Misener \protect\BIBand{} Floudas(2014)}]{misener2014antigone}
Misener R, Floudas CA (2014) {ANTIGONE}: algorithms for continuous/integer
  global optimization of nonlinear equations. \emph{Journal of Global
  Optimization} 59(2-3):503--526.

\bibitem[{Paneque et~al.(2018)Paneque, Gendron, Lurkin, Azadeh,
  \protect\BIBand{} Bierlaire}]{paneque2018lagrangian}
Paneque MP, Gendron B, Lurkin V, Azadeh SS, Bierlaire M (2018) A lagrangian
  relaxation technique for the demand-based benefit maximization problem.
  \emph{Proceedings of the 18th {Swiss} {Transport} {Research} {Conference}}
  (Ascona, Switzerland).

\bibitem[{Pardalos \protect\BIBand{} Schnitger(1988)}]{pardalos1988checking}
Pardalos PM, Schnitger G (1988) Checking local optimality in constrained
  quadratic programming is np-hard. \emph{Operations Research Letters}
  7(1):33--35.

\bibitem[{Pardalos \protect\BIBand{} Vavasis(1991)}]{pardalos1991quadratic}
Pardalos PM, Vavasis SA (1991) Quadratic programming with one negative
  eigenvalue is {NP}-hard. \emph{Journal of Global Optimization} 1(1):15--22.

\bibitem[{Sahinidis(1996)}]{sahinidis1996baron}
Sahinidis NV (1996) Baron: A general purpose global optimization software
  package. \emph{Journal of Global Optimization} 8(2):201--205.

\bibitem[{Schlicher \protect\BIBand{} Lurkin(2022)}]{schlicher2021stable}
Schlicher L, Lurkin V (2022) Stable allocations for choice-based collaborative
  price setting. \emph{European Journal of Operational Research} .

\bibitem[{Song \protect\BIBand{} Xue(2007)}]{song2007demand}
Song JS, Xue Z (2007) Demand management and inventory control for substitutable
  products. \textit{Working paper}.

\bibitem[{Soon(2011)}]{soon2011review}
Soon W (2011) A review of multi-product pricing models. \emph{Applied
  Mathematics and Computation} 217(21):8149--8165.

\bibitem[{Sumida et~al.(2019)Sumida, Gallego, Rusmevichientong, Topaloglu,
  \protect\BIBand{} Davis}]{sumida2019revenue}
Sumida M, Gallego G, Rusmevichientong P, Topaloglu H, Davis J (2019)
  Revenue-utility tradeoff in assortment optimization under the multinomial
  logit model with totally unimodular constraints. Technical report, Cornell
  University, Ithaca, NY.

\bibitem[{Train(2003)}]{train_2003}
Train KE (2003) \emph{Mixed Logit}, 138–154 (Cambridge University Press).

\bibitem[{van~de Geer \protect\BIBand{} den Boer(2022)}]{GeerBranch}
van~de Geer R, den Boer AV (2022) Price optimization under the finite-mixture
  logit model. \emph{Management Science} Ahead of Print,
  \urlprefix\url{http://dx.doi.org/10.1287/mnsc.2021.4272}.

\bibitem[{Zhang \protect\BIBand{} Lu(2013)}]{zhang2013assessing}
Zhang D, Lu Z (2013) Assessing the value of dynamic pricing in network revenue
  management. \emph{INFORMS Journal on Computing} 25(1):102--115.

\bibitem[{Zhang et~al.(2018)Zhang, Rusmevichientong, \protect\BIBand{}
  Topaloglu}]{zhang2018multiproduct}
Zhang H, Rusmevichientong P, Topaloglu H (2018) Multiproduct pricing under the
  generalized extreme value models with homogeneous price sensitivity
  parameters. \emph{Operations Research} 66(6):1559--1570.

\bibitem[{Zhen et~al.(2021)Zhen, de~Moor, \protect\BIBand{} den
  Hertog}]{zhen2021extension}
Zhen J, de~Moor D, den Hertog D (2021) An extension of the
  reformulation-linearization technique to nonlinear optimization.
  \emph{Available at Optimization Online} .

\bibitem[{Zhen et~al.(2022)Zhen, Marandi, de~Moor, den Hertog,
  \protect\BIBand{} Vandenberghe}]{zhen2018disjoint}
Zhen J, Marandi A, de~Moor D, den Hertog D, Vandenberghe L (2022) Disjoint
  bilinear optimization: A two-stage robust optimization perspective.
  \textit{Informs Journal on Computing}.

\end{thebibliography}
\newpage
\begin{appendices}
\section{Illustration of the objective functions of the continuous mixed logit model for the case study}\label{Appendix: figure for objective functions Mixed}
\begin{figure}[ht]
    \begin{subfigure}{\linewidth}
    \centering
    \includegraphics[width=13.5cm, height=6.5cm]{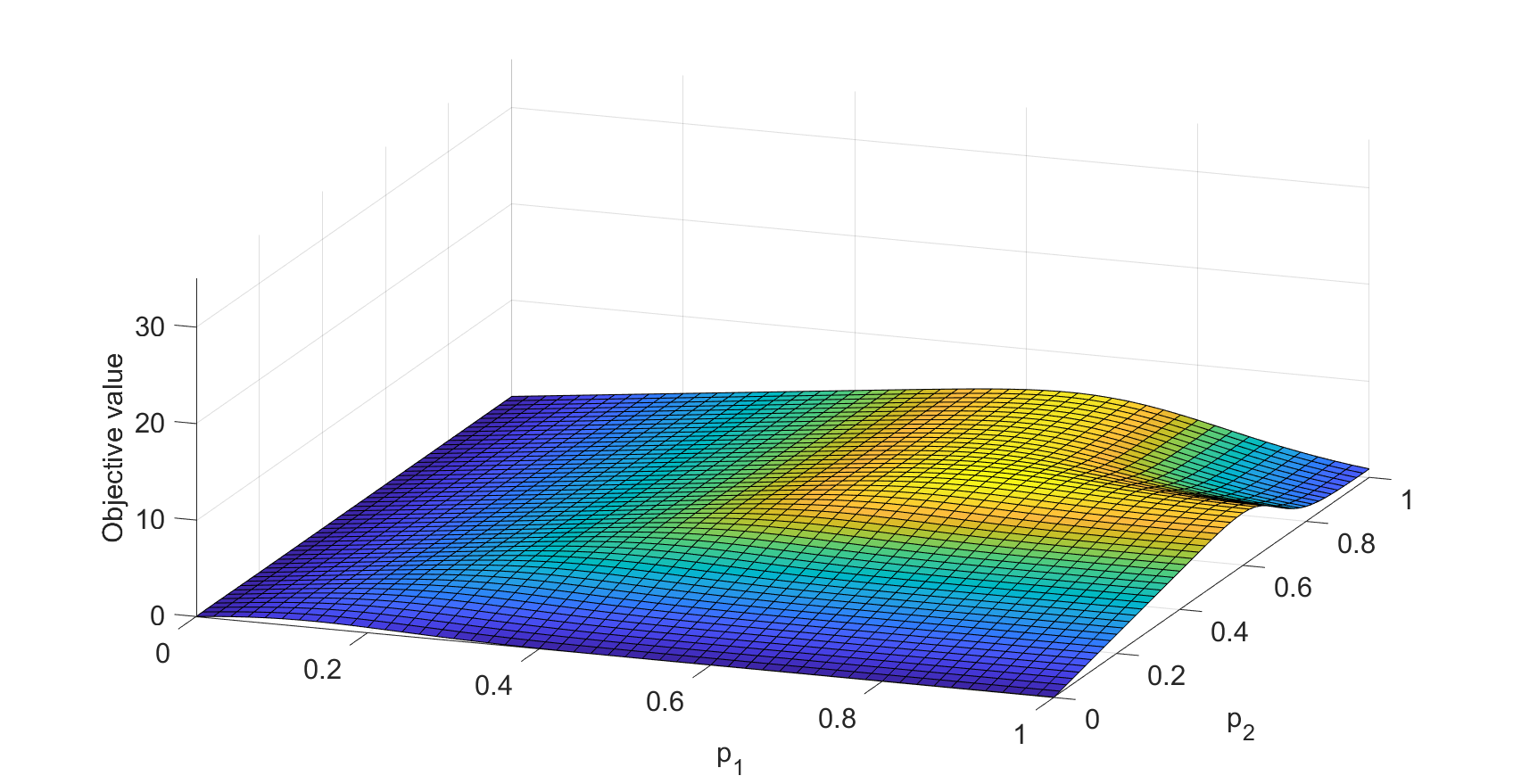}
    \subcaption{$N=10$}
    \label{fig:Mixedobjective n=10}
    \end{subfigure}
       \begin{subfigure}{\linewidth}
       \centering
    \includegraphics[width=13cm, height=6.5cm]{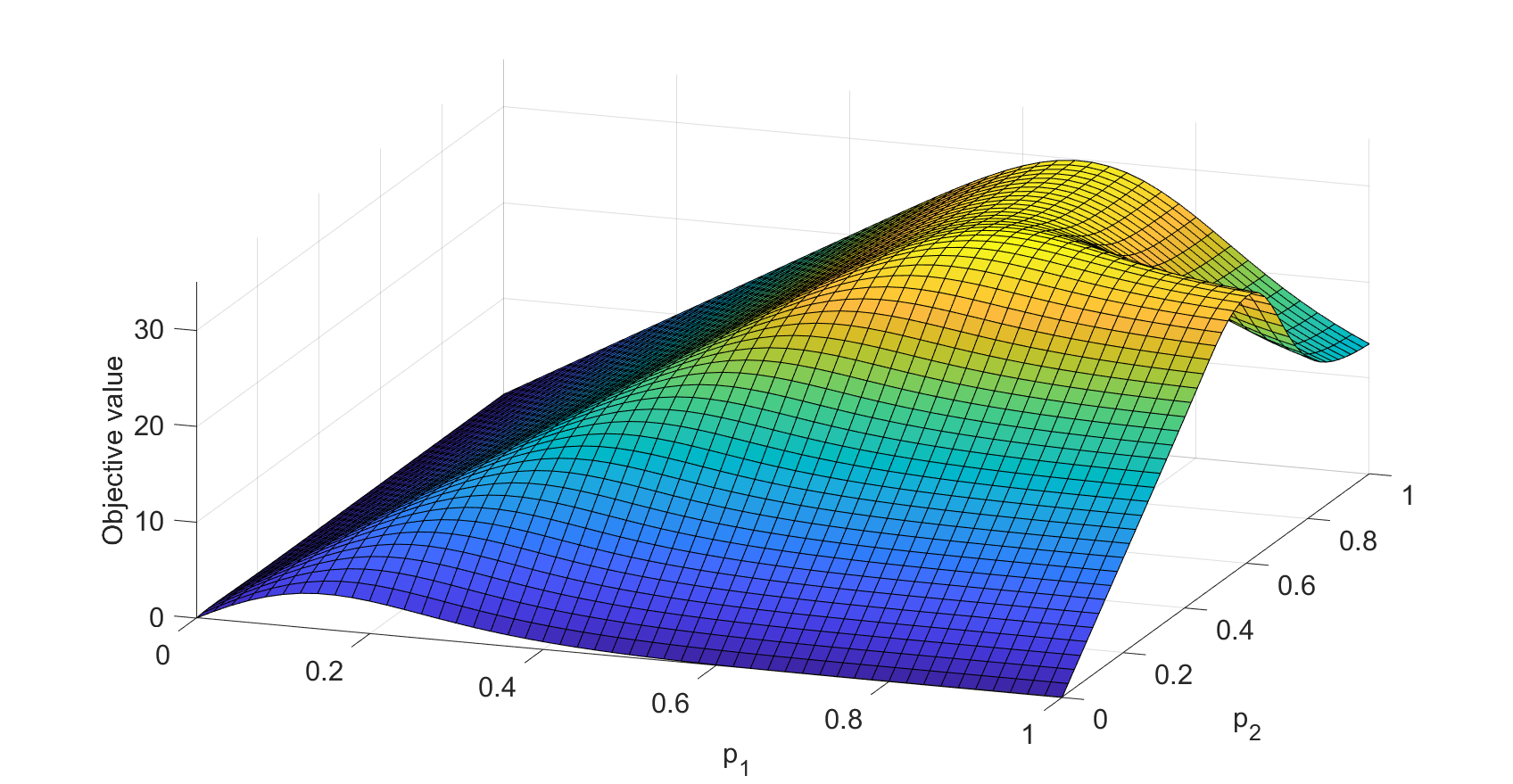}
    \subcaption{$N=50$}
    \label{fig:Mixedobjective n=50}
    \end{subfigure}
    \caption{Illustration of the objective function of \eqref{Eq:MixedMNL} for the parking choice model with $N=10$ and $N=50$ customers.}
    \label{fig:MixedLogit_10_objective}
\end{figure}
\end{appendices}
\end{document}